\documentclass[11pt]{article}
\usepackage{amscd,amssymb,stmaryrd}
\usepackage[fleqn]{amsmath}
\usepackage{graphicx}
\usepackage{subcaption}
\usepackage[utf8]{inputenc}
\usepackage[export]{adjustbox}
\usepackage{wrapfig}
\usepackage{amstext}
\usepackage{pstricks,pst-node,pst-plot,pst-coil}
\usepackage{amsthm}
\usepackage{algpseudocode}
\usepackage{algorithm}
\usepackage[title]{appendix}
\usepackage{amsmath}
\usepackage{color}
\usepackage{mathtools}
\usepackage{comment}
\usepackage{hyperref}
\usepackage[english]{babel}
\usepackage{booktabs}
\usepackage{mathrsfs}
\usepackage{tcolorbox}

\newcommand\norm[1]{\left\|#1\right\|}
\newcommand\abs[1]{\lvert#1\rvert}

\newcommand{\tand}{\quad \text{and} \quad}
\newtheorem{theorem}{Theorem}[section]

\newtheorem{example}[theorem]{Example}
\theoremstyle{remark}

\newcommand{\mca}{\mathcal{A}}
\newcommand{\mcr}{\mathcal{R}}

\newcommand{\cz}{\bar{u}}
\newcommand{\cd}{\tilde{u}}
\newcommand{\vz}{\bar{a}}
\newcommand{\vd}{\tilde{a}}

\newcommand{\veps}{\varepsilon}

\title{Nonlocal transport equations in multiscale media. Modeling, dememorization, and discretizations}
\author{Yalchin Efendiev\footnote{Department of Mathematics, Texas A\&M University, College Station, TX 77843, USA}, ~ 
Wing Tat Leung\thanks{Department of Mathematics, 
University of California Irvine, Irvine, CA 92697, USA}, ~ 
Wenyuan Li\footnote{Department of Mathematics, Texas A\&M University, College Station, TX 77843, USA}, ~ \\
Sai-Mang Pun\footnote{Department of Mathematics, Texas A\&M University, College Station, TX 77843, USA}, ~ and ~ Petr N. Vabishchevich\footnote{Nuclear Safety Institute, Russian Academy of Sciences, Moscow, Russia \& North-Caucasus Federal University, Stavrapol, Russia}}
\begin{document}
\maketitle
\begin{abstract}
In this paper, we consider a class of convection-diffusion equations
with memory effects. These equations arise as a result of homogenization 
or upscaling
of linear transport equations in heterogeneous media and play an important
role in many applications. First, we present a dememorization technique
for these equations. We show that the 
convection-diffusion equations
with memory effects can be written as a system of standard convection diffusion
reaction equations. This allows removing the memory term and simplifying
the computations. We consider a relation between dememorized equations and
micro-scale equations, which do not contain memory terms. We note that
dememorized equations differ from micro-scale equations and constitute
a macroscopic model. Next, we consider both implicit and partially explicit
methods. The latter is introduced for problems in multiscale media with
high-contrast properties.  Because of high-contrast, explicit methods
are restrictive and require time steps that are very small
(scales as the inverse of the contrast). We show that, by
 appropriately
decomposing the space, we can treat only a few degrees of freedom
implicitly and the remaining degrees of freedom explicitly. 
We present a stability analysis.
Numerical results are presented that confirm our theoretical
findings about partially explicit schemes applied to 
dememorized systems of equations.
\end{abstract}

\section{Introduction}
There are many problems that contain memory terms 
\cite{auriault1993deformable,christensen2012theory,pruss2013evolutionary,tartar1991memory}. One the well-known
example is macro-dispersion due to small scales and  reaction at
small scales. It has the form
\begin{eqnarray}
\label{eq:main1}
u_t + a(x) \cdot \nabla u = \int_0^t \nabla \cdot \left ( A(x,t,s) \nabla u (\tilde x(x,t,s),s) \right )ds.
\end{eqnarray}
Here, $u$ represents flow saturation in some porous medium,
the term $a(x)$ is a given velocity field at the macroscopic level, 
and $\tilde x = \tilde x (x,t,s)$ is a trajectory that depends on 
fine-scale heterogeneities. 
%The symbol $u_t$ is the derivative of $u = u(x,t)$ with respect to the time variable $t$ and the gradient operator $\nabla$ is related to the derivative with respect to the spatial variable $x$. 
The macro-dispersion (the term on the right hand side) 
is due to small-scale fluctuations of the velocity 
and the reaction at the
micro-scale. The equation \eqref{eq:main1} occurs in many porous
media related applications \cite{furman1997heat}, which include groundwater, petroleum engineering, 
and biomedical applications. The velocity and macro-dispersion terms are, in 
general, heterogeneous as the velocity fluctuations are upscaled over the
smallest scales. In this paper, our goal
is to show how to dememorize these types of problems and its relation to
homogenization, which does not contain memory related terms.

%Problems with memory. Describe. Relation to homogenization.
Solving \eqref{eq:main1} involves handling the memory terms
and saving all previous time information. This can be difficult 
especially for multiscale and nonlinear problems. There have been 
several approaches that dememorize the problems of a different form
\cite{vabishchevich2021approximate}. In this paper, 
we follow similar concepts and 
dememorize \eqref{eq:main1} and consider its relation to 
equations at the micro-scales, which do not contain memory terms.
In particular, we show that dememorized equations are, in some sense,
homogenized equations, not similar to equations at the micro-scale.
For example, the convection in dememorized problems 
occurs with averaged velocities. The diffusion is related 
to average quantity, which is represented by
$u$.

Dememorized equations constitute a system of coupled equations.
In particular, the main equation can be written as
$$
u_t + a\cdot \nabla u = v,
$$
where $v$ is due to perturbation from the average state. The equations
for $v$ are convection-diffusion-reaction types, where convection and
reaction effects occur with average rates.  This equation contains
a diffusion term, which depends on $u$. Note that the equations at 
the micro-scales are purely convection-reaction types. We discuss
the relation to microscale equations.

In this paper, we study dememorization and its discretization. We consider
two types of discretizations, namely implicit and partially explicit. The
latter is designed for multiscale problems based on a solution decomposition strategy \cite{efendiev2021temporal,efendiev2021splitting}. 
Because of the multiscale nature of the velocity and diffusion terms, 
one needs a very small time step
when performing explicit discretization. The time step depends 
on the contrast. In \cite{chung2021contrast}, a partially 
explicit approach was first proposed
for heterogeneous parabolic equations. The main idea of this approach 
is to handle some degrees of freedom implicitly, while the rest
explicitly. As a result, we identify a few degrees of freedom on a coarse
grid, that is much larger compared to spatial heterogeneities, 
and treat them implicitly. Our previous works (see 
also \cite{hu2021partially,li2021partially}) show that the
 resulting approach is stable
with appropriate decomposition of implicit and explicit components.
In particular, implicit components account for {\it fast} flows, while
explicit components account for {\it slow} flows. 
In this paper, we extend the partially explicit concept for the dememorized
equations of \eqref{eq:main1}. We use a spatial decomposition
of the solution following the constraint energy minimizing generalized
multiscale finite element method (CEM-GMsFEM) 
previously developed in \cite{chung2018constraint}. 
In this decomposition, the fast and slow components of the solution
are identified. Furthermore, we use implicit discretization
for fast components and explicit discretization for
slow components.

The rest of the paper is organized as follows. In Section \ref{sec:prelim}, we present some preliminaries for the model problem. In particular, we dememorize and derive the coupled system equivalent to the nonlocal equation with memory effects. In Section \ref{sec:discretize}, we derive the numerical discretization schemes for the model problem. Numerical experiments are presented in Section \ref{sec:numerics}. Concluding remarks are drawn in Section \ref{sec:conclusion}. 
In the Appendices, we present some remarks related to homogenization,
constructing  multiscale spaces, and
stability estimates.

\section{Preliminaries} \label{sec:prelim}
In this section, we present some preliminaries for the nonlocal transport equations
that arise in porous media applications. In particular, we introduce the model problem considered in this work, and derive the coupled system using the technique of dememorization. 
\subsection{Model Problem}
Let $\Omega \subset \mathbb{R}^d$ ($d \in \{ 2, 3\}$) be a bounded domain and $T>0$ be a given terminal time. 
We consider the following boundary value problem  with memory effects: find $u(x,t)$ such that 
\begin{align}
%\begin{split}
 u_t (x,t)+ a(x) \cdot \nabla u(x,t) &=  \int_0^t \nabla \cdot \left ( A(x,t,s) \nabla u (\tilde x(x,t,s),s) \right )ds \text{in}~ \Omega \times (0,T], \label{eqn:model} \\
 u(x,0) &= u_0(x) \text{in} ~ \Omega, \nonumber\\
 u(x,t) &=0\  \text{on} \ \partial\Omega, \nonumber
 %\end{split}
\end{align}
where $a : \mathbb{R}^d \to \mathbb{R}^d$ is a vector-valued function, 
$A(x,t,s)$ represents a macro-dispersion coefficient,  %$f \in L^2(0,T;L^2(\Omega))$ is a source term, 
and $\tilde x(x,t,s)$ is a trajectory that satisfies
\begin{eqnarray}
{d\over ds} \tilde x(x,t,s) = \widetilde{a}(\tilde x), \quad  \tilde x (x,t,t) = x.
\label{eqn:trajectory}
\end{eqnarray}
That is, the trajectory with the velocity $\widetilde{a}(\tilde x)$ is such that at the time $t$ it reaches the point $x$. %\textbf{is $\tilde{x}$ has a variable $x$?}   
%See Appendix \ref{app:perturb} for the discussion of the upscaling procedure. 
The trajectory has the following explicit expression $\tilde x(x,t,s) = x- (t-s) \widetilde a $ when  $\widetilde a $ is a constant.
The kernel $A(x,t,s)$ is assumed to be in terms of
the exponential term due to the reaction at the micro-scale. 
In this case, we assume that 
\begin{eqnarray}
A(x,t,s) = \kappa (x) e^{-\beta(t-s)}
\label{eqn:kernel-form}
\end{eqnarray}
for some permeability tensor $\kappa(x)$. Such reaction kernel can be derived from homogenization (see Appendix \ref{app:hom}).
More generally, the kernel $A(x,t,s)$ has the form 
\begin{eqnarray}
A(x,t,s) = \sum_{i=1}^M \kappa_i(x) e^{-\beta_i(t-s)},
\label{eqn:kernel-form-2}
\end{eqnarray}
where $M$ is a positive integer. In this case, the functions $\kappa_i$'s are some heterogeneous fields. 

The model problem \eqref{eqn:model} is well-known 
in the sense that it can be derived from upscaling
(see Appendix  \ref{app:perturb})
of some micro-scale transport equations containing heterogeneous 
 velocity fields \cite{dagan1984solute,zhang1995asymptotic}. 
In this case, the solution of the macroscopic equation is an 
average of the microscopic solution, where microscopic equations
do not contain memory terms. In this work, we show that one can 
re-write the macroscopic equations without memory. However, the resulting
macroscopic equation is different from microscale equations (without memory).
The form of macroscopic diffusion $A(x,t,s)$ is similar to the one
obtained from upscaling.

\subsection{Dememorization} \label{sec:dememorize}
In this section, we apply the technique of dememorization for the problem \eqref{eqn:model} with the kernel function having the form \eqref{eqn:kernel-form-2} and the trajectory \eqref{eqn:trajectory}. In particular, we introduce auxiliary variables $\{v_i\}_{i=1}^M$ and derive the coupled system for the main variable $u$ and the auxiliary variables. % $\{v_i\}_{i=1}^M$. 
%In the following, we assume that the kernel $A(x,t,s)$ has the form of \eqref{eqn:kernel-form-2}. 
%The case with a more general kernel in \eqref{eqn:kernel-form-2} can also be derived in a straightforward manner by simply introducing more auxiliary variables. 

The dememorization starts with the following auxiliary variables. For $i\in \{ 1,2,\cdots,M\}$, we define 
\begin{eqnarray}
\label{eq:u-rel-v}
v_{i} (x,t):= \int_0^t e^{-\beta_i(t-s)} u (\tilde x(x,t,s),s)ds \quad \text{for any} ~ (x,t) \in \Omega \times (0,T].
\end{eqnarray}
Note that, from the original equation \eqref{eqn:model}, we have 
$$u_t + a(x) \cdot \nabla u = \sum_{i=1}^M \nabla \cdot (\kappa_i (x) \nabla v_{i}).$$
On the other hand, 
taking derivative (with respect to $t$) of $v_i$, we obtain 
$$\left(v_{i}\right)_t + \beta_i v_{i} + \widetilde{a} \cdot \nabla v_{i} =u.$$
Therefore, we obtain the following coupled system without memory effects: 
\begin{eqnarray}
\begin{split}
\left (v_{i}\right)_t + \beta_i v_{i} + \widetilde{a} \cdot \nabla v_{i} &=u& \quad \text{in} ~ \Omega \times (0,T],\; i \in \{ 1,2,\cdots,M \},\\
u_t + a \cdot \nabla u-\sum_{i=1}^M \nabla \cdot (\kappa_i \nabla v_{i}) &= 0 & \quad \text{in} ~ \Omega \times (0,T], \\
u(x,0) &= u_0(x) &\quad \text{in} ~ \Omega, \\
u(x,t) &= 0 &\quad \text{in} ~ \partial \Omega.
\label{eqn:model-coupled}
\end{split}
\end{eqnarray}
The boundary condition for each $v_i(x,t)$ is defined via $u$ using 
\eqref{eq:u-rel-v} while assuming $u$ outside $\Omega$ is zero. 
We remark that if the permeability tensor has the form $\kappa = \kappa_{11}\left(a \otimes a \right )$ for some bounded function $\kappa_{11}$ and the conditions $\beta \kappa_{11}+ \nabla \kappa_{11} \cdot a \geq 0$ and $a \cdot \mathbf{n}_{\partial \Omega} = 0$ hold, where $\mathbf{n}_{\partial \Omega}$ is the unit outward normal vector to the boundary $\partial \Omega$, then one can show that the continuous problem \eqref{eqn:model-coupled} is stable. 
See Appendix \ref{app:stable} for more details on the stability of the continuous problem \eqref{eqn:model-coupled}. 
The stability analysis for the discretized convection-diffusion model with memory effects is challenging and will be one of our future works. In the following, we develop numerical discretization scheme for \eqref{eqn:model-coupled} and provide a stability estimate for the case when $a = \widetilde{a} \equiv \mathbf{0}$.

For the numerical discretization, we introduce the variational 
formulation of the problem \eqref{eqn:model-coupled}. To this aim, we define $\Gamma \subset \partial \Omega$ such that $\Gamma := \{ x \in \partial \Omega: \widetilde a(x) \cdot \mathbf{n}_{\partial \Omega} (x) <0\}$ and we write $H_\Gamma^1 (\Omega) := \{ v \in H^1(\Omega) : v \vert_{\Gamma} = 0 \}$. We assume that $\Gamma$ has positive measure. 
The corresponding variational problem reads as follows: 
Find $u(\cdot, t) \in H_0^1(\Omega)$ and 
$\{v_i(\cdot,t)\}_{i=1}^M \subset  H_{\Gamma}^1(\Omega)$ 
 such that 
\begin{eqnarray}
\begin{split}
(\left (v_{i} \right )_t,\phi) +(\beta_i v_{i},\phi) + (\widetilde{a} \cdot \nabla v_{i}, \phi) &=(u,\phi) &\quad \forall \phi \in H_{\Gamma}^1(\Omega),\quad i \in \{ 1,2,\cdots,M \},\\
(u_t, \psi)+(a \cdot \nabla u,\psi) + \sum_{i=1}^M \mca_i(v_{i},\psi) &=0 &\quad \forall \psi \in H_0^1(\Omega).  
\end{split}
\label{eqn:var-model}
\end{eqnarray}
We denote $(\cdot,\cdot)$ the inner product in $L^2(\Omega)$ and $\mca_i(\cdot,\cdot)$ is defined to be 
$$ \mca_i(\phi,\psi) := \int_\Omega \kappa_i \nabla \phi \cdot \nabla \psi ~ dx $$
for any $\phi \in H^1(\Omega)$ and $\psi \in H^1(\Omega)$. We denote $\norm{\cdot}$ the $L^2$ norm induced by the inner product $(\cdot,\cdot)$ and we write $\norm{\cdot}_{\mca_i} := \sqrt{\mca_i(\cdot,\cdot)}$. 
When $M = 1$, we simply write 
$$ \mca(\phi,\psi) := \int_\Omega \kappa_1 \nabla \phi \cdot \nabla \psi ~ dx $$
for any $\phi \in H^1(\Omega)$ and $\psi \in H^1(\Omega)$ and the corresponding energy norm is written as $\norm{\cdot}_{\mca} := \sqrt{\mca (\cdot,\cdot)}$. 
For $i\in \{ 1,2,\cdots,M \}$, we assume that $\beta_i$ is a positive constant and $\kappa_i \in L^\infty(\Omega; \mathbb{R}^{d\times d})$ is a permeability tensor that fulfils $0< C_0 \leq \xi^T \kappa_i \xi \leq C_{\infty} < \infty$ for any $\xi \in \mathbb{R}^d$ with $\abs{\xi} = 1$ (with $\abs{\cdot}$ being the usual Euclidean norm in $\mathbb{R}^d$). 
We may define the operator $\mathcal{R}_i: H_0^1(\Omega) \to H^{-1} (\Omega)$ such that 
$$(\mathcal{R}_i \phi ,\psi) := \mca_i(\phi,\psi)$$
for any $\phi \in H^1(\Omega)$ and $\psi \in H^1(\Omega)$. 

Next, we derive an energy estimate for the solution in the absence 
of convection (see Appendix \ref{app:stable} for the case with convection).
We remark that assuming $a = \widetilde{a} \equiv \mathbf{0}$, if we take $\psi = u$ in \eqref{eqn:var-model} and integrate over $(0,T]$, we obtain 
\begin{eqnarray*}
\begin{split}
0 & =  \int_{0}^{T} \left [ (u_t,u) + \sum_{i=1}^M (\mathcal{R}_iv,u) \right ]dt \\
& = \int_{0}^{T}\left [ (\partial_{t}u,u) + \sum_{i=1}^M (\mathcal{R}_i v_{i}, \beta v_{i}) + \sum_{i=1}^M (\mathcal{R}_i v_{i}, \partial_t v_{i})\right ]dt\\
& =  \frac{1}{2} \int_0^T \frac{d}{dt} \left ( \norm{u}^2 + \sum_{i=1}^M \norm{v_{i}}_{\mca_i}^2  \right ) dt+ \sum_{i=1}^M \int_0^T \mca(v_{i}, \beta v_{i})dt\\
& = \cfrac{1}{2}\left (\|u(\cdot,T)\|^{2}+\sum_{i=1}^M\norm{v_{i}(\cdot,T)}_{\mca_i}^{2}-\|u_0\|^{2} - \sum_{i=1}^M\norm{v_{i}(\cdot,0)}_{\mca_i}^2 \right ) + \sum_{i=1}^M \beta_i \int_0^T \norm{v_{i}(\cdot,t)}_{\mca_i}^2  dt.
\end{split}
\end{eqnarray*}
Therefore, we have 
$$
E(u,v;M;T) := \|u(\cdot,T)\|^{2}+\sum_{i=1}^M\|v_{i}(\cdot,T)\|_{\mca_i}^{2} \leq \|u_0\|^{2}+\sum_{i=1}^M\|v_i(\cdot,0)\|_{\mca_i}^{2} =: E(u,v;M;0)
$$
for the function $E(u,v;M;t) := \norm{u(\cdot,t)}^2 + \sum_{i=1}^M \norm{v_{i}(\cdot,t)}_{\mca_i}^2$. It implies that the continuous problem \eqref{eqn:var-model} (when $a = \widetilde{a} \equiv \mathbf{0}$) is stable with respect to the function $E(\cdot,\cdot;M;t)$. 

\section{Numerical Discretizations} \label{sec:discretize}
In this section, we set $M=1$ and  develop the numerical discretization for \eqref{eqn:var-model} and provide stability estimate for the numerical schemes. We present the temporal discretizations using the implicit scheme and the recently developed partially explicit scheme based on a space decomposition strategy. 
We assume that some finite dimensional spaces $V_H \subset H_0^1(\Omega)$ and $W_H \subset H_{\Gamma}^1(\Omega)$ based on some (coarse-grid) partition for the domain are developed and we perform spatial discretization using Galerkin method with the ansatz spaces $V_H$ and $W_H$. The coarse-grid spaces $V_H$ and $W_H$ are defined via the recently developed CEM-GMsFEM for multiscale problems. See Appendix \ref{SpaceConstructSection} for more details of the construction and definition. 

\subsection{Semi-implicit Scheme}
In this section, we develop the implicit-in-time fully discretization for the problem \eqref{eqn:var-model}.
%\marginpar{This section needs more careful writing via bilinear forms also!}
To this aim, we introduce a temporal partition $\{ t^n \}_{n=1}^{N_T}$ with $t^n = n\Delta t$ ($k \in \{ 0, 1, \cdots, N_T \}$) and $T = N_T \Delta t$; we also denote $v^n = v(\cdot, t^n)$ and $u^n = u(\cdot,t^n)$ for any $n \in \{ 0,1,\cdots, N_T \}$. 
We remark that the convection term is computed explicitly in our discretization. 
%The implicit scheme for (\ref{eqs1}) is 
%\begin{eqnarray*}
%    \cfrac{v^{n+1}-v^n}{\Delta t} + \beta v^{n} + \widetilde{a} \cdot \nabla v^n - u^{n+1} = 0, \\
%    \cfrac{u^{n+1}-u^n}{\Delta t} + a \cdot \nabla u^n - \nabla \cdot (\kappa (x) \nabla v^{n+1}) = 0.
%\end{eqnarray*}
%For the spatial discretization, we apply the finite element method with some finite element space $V_H \subset V$ that is based on a coarse grid of the spatial domain. 
The (implicit) fully discretization reads as follows: find $\{ u_H^k \}_{k=1}^{N_T} \subset V_H$ and $\{ v_H^k \}_{k=1}^{N_T}  \subset W_H$ such that the following system holds 
\begin{eqnarray}
\begin{split}
\left (\cfrac{v_H^{n+1}-v_H^n}{\Delta t}, \phi \right ) + \beta (v_H^{n}, \phi ) + (\widetilde{a} \cdot \nabla v_H^n, \phi ) - (u_H^{n+1}, \phi) & = 0 & \quad \forall \phi \in W_H, \\
\left (\cfrac{u_H^{n+1}-u_H^n}{\Delta t}, \psi \right ) + (a \cdot \nabla u_H^n,\psi) + \mca(v_H^{n+1}, \psi) &= 0 & \quad \forall \psi \in V_H.
\end{split}
\label{eqn:fully-implicit}
\end{eqnarray}
for any $n \in \{0, 1,\cdots, N_T -1\}$. 
The terms $u_H^0$ and $v_H^0$ are obtained from the initial conditions in the sense that 
$$(v_H^0, \phi) = (v_0,\phi) \quad \forall \phi \in W_H \tand (u_H^0, \psi) = (u_0, \psi) \quad \forall \psi \in V_H.$$
We remark that assuming $a = \widetilde{a} \equiv \mathbf{0}$, one can show the stability of the fully implicit scheme \eqref{eqn:fully-implicit}. 
Let $\psi = u_H^{n+1}$ in \eqref{eqn:fully-implicit} and we have 
\begin{eqnarray*}
\begin{split}
    0 &= \cfrac{1}{\Delta t}(u_H^{n+1}-u_H^n,u_H^{n+1}) +  (\mcr v_H^{n+1},u_H^{n+1}) \\ 
    &= \cfrac{1}{2\Delta t}(\|u_H^{n+1}\|^2-\|u_H^n\|^2+\|u_H^{n+1}-u_H^n\|^2 ) + \frac{1}{\Delta t} (\mcr v_H^{n+1} , v_H^{n+1} - v_H^n) + \beta (v_H^{n+1}, v_H^{n+1})   \\
%    &\geq \cfrac{1}{2\Delta t} (\|u_H^{n+1}\|^2-\|u_H^n\|^2) + \tilde{b}^{-1} a(v_H^{n+1},v_H^{n+1}) - \tilde{b}^{-1} e^{-b dt} a(v_H^{n+1}, v_H^n)  \\
%    &\geq \cfrac{1}{2\Delta t} (\|u_H^{n+1}\|^2-\|u_H^n\|^2) + \tilde{b}^{-1} (1-e^{-b \Delta t}) \norm{v_H^{n+1}}_a^2 + \tilde{b}^{-1} e^{-b\Delta t} a(v_H^{n+1}, v_H^{n+1}-v_H^n)  \\ 
    &\geq \cfrac{1}{2\Delta t} (\|u_H^{n+1}\|^2-\|u_H^n\|^2) + \beta \|v_H^{n+1}\|^2 + \cfrac{1}{2\Delta t} (\| v_H^{n+1}\|_{\mca}^2 - \|v_H^n\|_{\mca}^2 + \|v_H^{n+1} - v_H^n\|_{\mca}^2)  \\ 
    &\geq  \cfrac{1}{2\Delta t} \left (\|u_H^{n+1}\|^2-\|u_H^n\|^2 \right ) +  \cfrac{1}{2\Delta t}  \left ( \|v_H^{n+1}\|_{\mca}^2 - \|v_H^n\|_{\mca}^2 \right ).
\end{split}
\end{eqnarray*}
Thus, we have shown that, for any $n \in \{1,\cdots, N_T\}$, 
\begin{eqnarray*}
E^n(u_H,v_H) \leq E^0(u_H, v_H), \quad \text{where} ~ E^n (u_H, v_H) := \|u_H^{n}\|^2 + \|v_H^{n}\|_{\mca}^2 .%\leq \|u_H^0\|^2 + \tilde{b}^{-1} e^{-b\Delta t} \Delta t \|v_H^0\|_a^2.
\end{eqnarray*}
This shows the stability for the case with  $a = \widetilde{a} \equiv \mathbf{0}$. In Appendix \ref{app:stable}, we give a stability proof for 
more general case.

\subsection{Partially Explicit Splitting Scheme} 
In this section, we first briefly review the recently developed partially explicit splitting scheme and apply this scheme for discretizing \eqref{eqn:var-model}. The partially explicit splitting scheme is based on a solution decomposition strategy for the coarse spaces $V_H$ and $W_H$. We assume each ansatz space can be written as a direct sum of two subspaces: 
$V_H = V_H^1 \oplus V_H^2$ and $W_H = W_H^1 \oplus W_H^2$; we seek approximations in these ansatz spaces. In particular, the component in the first subspace $V_H^1$ (resp. $W_H^1$) will be treated implicitly during the evolution while the components in the second subspace $V_H^2$ (resp. $W_H^2$) will be computed in an explicit manner. An enhancement in terms of computational efficiency can be achieved within this setting of implicit-explicit formulation. 
%In the following, the subscript $i \in \{1,2\}$ of the function $v_{H,i}$ indicates that $v_{H,i} \in V_H^i$ to avoid any confusion to the auxiliary variables $v_{i}$ in the dememorization. 
%The precise construction of this ansatz spaces $V_H^1$ and $V_H^2$ is presented in Appendix. 

%We use the finite element method to solve the above two schemes. 
%Let $W_H$ be the finite element space and $W_H = W_{H,1} \oplus W_{H,2}$. 
With these ansatz spaces and the specific subspace decomposition, we can write $u_H^n = u_{H,1}^n + u_{H,2}^n$ and $v_{H}^n = v_{H,1}^n + v_{H,2}^n$ for any $n \in \{ 0, 1,\cdots, N_T \}$. 
%Likewise, we have $u_H = u_{H,1} + u_{H,2}$ and $v_H = v_{H,1} + v_{H,2}$ with $u_{H , i}, v_{H,i}\in W_{H,i}$. 
%In the partially explicit scheme, we split the space $W$ into two subspaces, i.e. $W = W_1 \oplus W_2$. Then, we have $v = v_1 + v_2$ and $u = u_1 + u_2$ where $v_i, u_i \in W_i \;(i=1,2)$.
The partially explicit splitting scheme reads as follows: find $\{ u_{H,i}^n \}_{n=1}^{N_T} \subset V_H^i$ and $\{ v_{H,i}^n \}_{n=1}^{N_T} \subset W_H^i$ for $i \in \{ 1, 2 \}$ such that 
the following system holds 
\begin{comment}
\begin{eqnarray*}
    \cfrac{v_1^{n+1}-v_1^n}{\Delta t} + \beta v_1^{n} + \lambda \cdot \nabla v_1^n - u_1^{n+1} = 0, \\
    \cfrac{v_2^{n+1}-v_2^n}{\Delta t} + \beta v_2^{n} + \lambda \cdot \nabla v_2^n - u_2^{n} = 0, \\
    \cfrac{u_1^{n+1}-u_1^n}{\Delta t} + \cfrac{u_2^{n+1}-u_2^{n}}{\Delta t} + a \cdot \nabla u^n - \nabla \cdot (\kappa (x) \nabla v^{n+1}) = 0,\\
    \cfrac{u_1^{n+1}-u_1^n}{\Delta t} + \cfrac{u_2^{n+1}-u_2^{n}}{\Delta t} + a \cdot \nabla u^n - \nabla \cdot (\kappa (x) \nabla v^{n+1}) = 0.
\end{eqnarray*}
\end{comment}
%Let $(\cdot,\cdot)$ be the $L^2$ inner product and $B(v_H,w_H) = \int_{\Omega} \kappa(x) \nabla v_H \cdot \nabla w_H$. 
%The partially explicit scheme is finding $v_H\in W_H$, $u_{H,1} \in W_{H,1}$ and $u_{H,2} \in W_{H,2}$ such that
\begin{eqnarray}
\begin{split}
\left(\cfrac{v_{H,1}^{n+1}-v_{H,1}^n}{\Delta t}, \phi_1 \right) + \beta (v_{H,1}^{n},\phi_1) + (\widetilde{a} \cdot \nabla v_{H,1}^n,\phi_1) - (u_{H,1}^{n+1},\phi_1) &= 0 & \quad \forall \phi_1 \in W_{H}^{1}, \\
\left(\cfrac{v_{H,2}^{n+1}-v_{H,2}^n}{\Delta t}, \phi_2 \right) + \beta (v_{H,2}^{n},\phi_2) + (\widetilde{a} \cdot \nabla v_{H,2}^n,\phi_2) - (u_{H,2}^{n},\phi_2)  & = 0& \quad \forall \phi_2\in W_{H}^{2}, \\
\left(\cfrac{u_{H,1}^{n+1}-u_{H,1}^n}{\Delta t} + \cfrac{u_{H,2}^{n+1}-u_{H,2}^{n}}{\Delta t} ,\psi_1\right) + (a \cdot \nabla u_{H}^n, \psi_1) + \mca(v_{H}^{n+1}, \psi_1) &= 0& \quad \forall \psi_1 \in V_{H}^{1},\\
\left(\cfrac{u_{H,1}^{n+1}-u_{H,1}^n}{\Delta t} + \cfrac{u_{H,2}^{n+1}-u_{H,2}^{n}}{\Delta t}, \psi_2\right) + (a \cdot \nabla u_{H}^n, \psi_2) + \mca(v_{H}^{n+1}, \psi_2) &= 0& \quad \forall \psi_2 \in V_{H}^{2},
\end{split}
\label{eqn:part-exp}
\end{eqnarray}
for any $n \in \{0, 1,\cdots, N_T -1\}$. 
For the case of pure reaction (i.e., $a = \widetilde{a}=\mathbf{0}$), we can derive a stability estimate for the above-mentioned partially explicit splitting scheme \eqref{eqn:part-exp}. To this aim, we define a constant $\gamma \in (0,1)$ such that 
\begin{eqnarray}
\gamma:=\sup_{v_{1}\in V_H^{1},~v_{2}\in V_H^{2}}\cfrac{(v_{1},v_{2})}{\|v_{1}\|\|v_{2}\|}.
\label{eqn:gamma}
\end{eqnarray}
For the pure reaction case, the partially explicit splitting scheme \eqref{eqn:part-exp} is stable under appropriate assumptions on the subspaces $V_{H}^1$ and $V_H^2$. The stability estimate for the general convection-diffusion case is left as future work. 
\begin{theorem}[Stability estimate of pure reaction case]
\label{thm:pe}
Assume that $a = \widetilde{a} \equiv \mathbf{0}$. Let $\gamma$ be defined in \eqref{eqn:gamma}. Suppose that the temporal step size $\Delta t$ satisfies 
\begin{eqnarray}
\Delta t\leq \beta(1-\gamma) \inf_{v \in V_H^2} \frac{\norm{v}^2}{\norm{v}_{\mca}^2}.
\label{eqn:pe-time-step}
\end{eqnarray}
%$$\cfrac{(1-\gamma)}{\Delta t} \geq \cfrac{1}{b}\sup\limits_{v\in V_H^2} \cfrac{\|v\|_{a}^{2}}{\|v\|^2}.$$
Then, the solutions $u_H^n = u_{H,1}^n + u_{H,2}^n$ and $v_{H}^n = v_{H,1}^n + v_{H,2}^n$ obtained from \eqref{eqn:part-exp} satisfy the following stability estimate
%\begin{align*}
%&\cfrac{e^{-bdt}}{2\tilde{b}}\|v^{n+1}\|_{a}^{2} + \cfrac{1}{2dt} \|u^{n+1}\|^2\leq \cfrac{e^{-bdt}}{2\tilde{b}}\|v^{n}\|_{a}^{2} + \cfrac{1}{2dt} \|u^n\|^2 
%\end{align*}
$$ \tilde E^n(u_H, v_H) \leq \tilde E^0(u_H,v_H)$$
for any $n \in \{ 1,\cdots, N_T \}$, where $\tilde E^n (u_H,v_H) := \displaystyle{\norm{u_H^n}^2 + \sum_{i=1}^2 \norm{v_{H,i}^{n}}_{\mca}^2}$ is the discrete energy function. 
\end{theorem}
The proof of this result is given in Appendix \ref{app:proof}.

\section{Numerical Experiments} \label{sec:numerics}
In this section, we perform some numerical experiments using the discretization schemes discussed in the previous section. 
In all the experiments, we set the spatial domain to be $\Omega = (0,1)^2$. We set the velocity fields to be $\widetilde{a} = (0.05,0)^T$ and $a = (0.1,0)^T$. 
The spatial domain is partitioned into uniform square elements with mesh size $H = \sqrt{2}/10$ to form a coarse grid. Next, for each coarse element from the coarse partition, we further divide it into $10 \times 10$ uniform square elements so that the mesh size $h$ of the fine grid is $h = \sqrt{2}/100$. 
We equip the variable $u$ with the homogeneous Neumann boundary condition on the whole boundary $\partial \Omega$. 

%The fine mesh size is $\cfrac{\sqrt{2}}{100}$ while the coarse mesh size is $\cfrac{\sqrt{2}}{10}$. We let the final time $T=0.2$. $\kappa$ is the permeability field that we use. In all numerical examples, we apply the Neumann boundary condition, \[ \nabla u \cdot \vec{n} = 0,\quad \text{with $\vec{n}$ being the outward normal vector.} \]
For the ansatz space, we choose three local auxiliary functions (i.e., $L_i = J_i = 3$ for each $i\in \{ 1, \cdots, N_e \}$ with $N_e = 100$; see Appendix \ref{SpaceConstructSection} for more details) to form the local auxiliary space in each coarse element so that the dimensions of $V_{H}^{1}$ and $V_H^2$ are $300$. 
We take the oversampling parameter to be $m = 4$. 
Based on the fine grid and implicit temporal discretization, we compute a numerical approximation which serves as a reference solution. 
In the following, we compute three different numerical approximations and compare them with the reference solution in terms of $L^2$ error: 
\begin{enumerate}
\item The first approximation is obtained using only the first ansatz space $V_H^1$ with the implicit temporal discretization (i.e., solving \eqref{eqn:fully-implicit} with the space $V_H^1$). We refer to this approximation as {\it implicit CEM}. 

\item Combining the additional ansatz space $V_H^2$, we compute the second approximation over the space $V_H = V_H^1 \oplus V_H^2$ via implicit scheme (i.e., solving \eqref{eqn:fully-implicit} with $V_H$). We refer to this approximation as {\it implicit CEM with additional bases}. 

\item The third one is computed by solving the partially explicit splitting scheme \eqref{eqn:part-exp}. We refer to this approximation as {\it partially explicit splitting CEM}. 
\end{enumerate}
%The implicit fine grid solution is calculated using the fine grid basis functions and via the implicit scheme. The implicit CEM solution is calculated implicitly using CEM bases. The implicit CEM solution with additional bases is obtained by the implicit scheme with basis functions from both $V_{H}^1$ and $V_{H}^2$. The partially explicit solution is calculated by the partially explicit scheme proposed above. We use the implicit fine grid solution as our reference solution and calculate the relative $L^2$ error for the other three schemes. 
From these numerical examples, we find that the partially explicit scheme can achieve similar accuracy as the fully implicit scheme with less computing cost at each time level. 

\begin{example} \label{exp:1}
In the first example, we set the initial condition to be 
$u_0(x) = \sin(\pi x_1) \sin(\pi x_2)$ for any $x = (x_1, x_2) \in \Omega$. %See Figure \ref{NRfigic1} (right) for a graphical demonstration for the initial condition. 
Let $T = 0.05$ and the temporal step size is $\Delta t = T/100 = 5 \times 10^{-4}$. 
The permeability field $\kappa$ used in this example is depicted in Figure \ref{NRfigic1}. % (left). 

In Figure \ref{NRfigic2}, we present the profiles of the three types of solutions at the terminal time - the reference solution (the implicit fine grid solution), the implicit CEM solution with additional bases, and the partially explicit splitting CEM solution. 
The relative $L^2$ error against time is presented in Figure \ref{NRfigic3}. 
Despite the differences among these profiles of the numerical approximations, 
%Although the profiles of the implicit CEM solution with additional bases and the partially explicit solution are a little bit different from the profile of the reference solution. 
the relative $L^2$ error is around $8\%$ using the partially explicit splitting scheme, which is relatively small and acceptable. 
Besides, from Figure \ref{NRfigic3}, the error curves for the implicit CEM solution with additional bases and the partially explicit solution nearly coincide. This implies that one can achieve the same level of accuracy using the proposed partially explicit splitting scheme as the implicit CEM scheme with additional basis functions. 

\begin{figure}[H]
\centering

\includegraphics[width = 8cm]{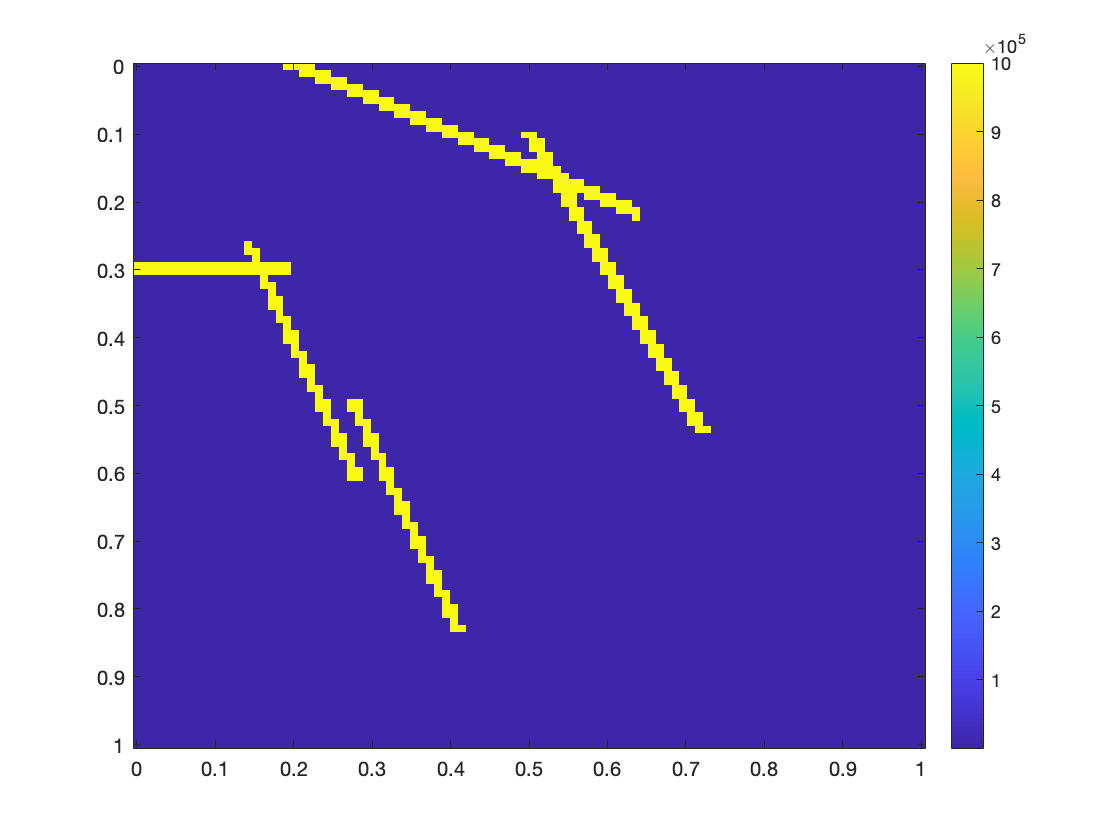}
\caption{Permeability in Example \ref{exp:1}.}
\label{NRfigic1}
\end{figure}

\begin{figure}[H]
\centering
\includegraphics[width = 5.2cm]{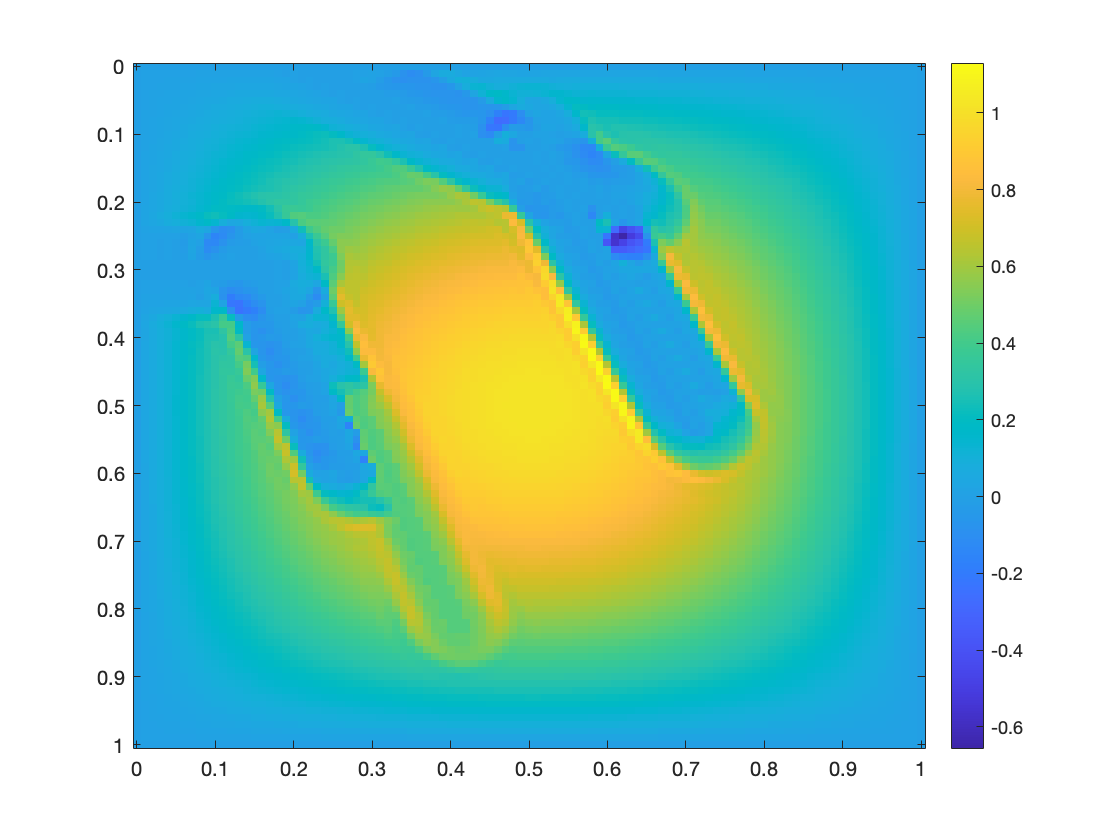}
\includegraphics[width = 5.2cm]{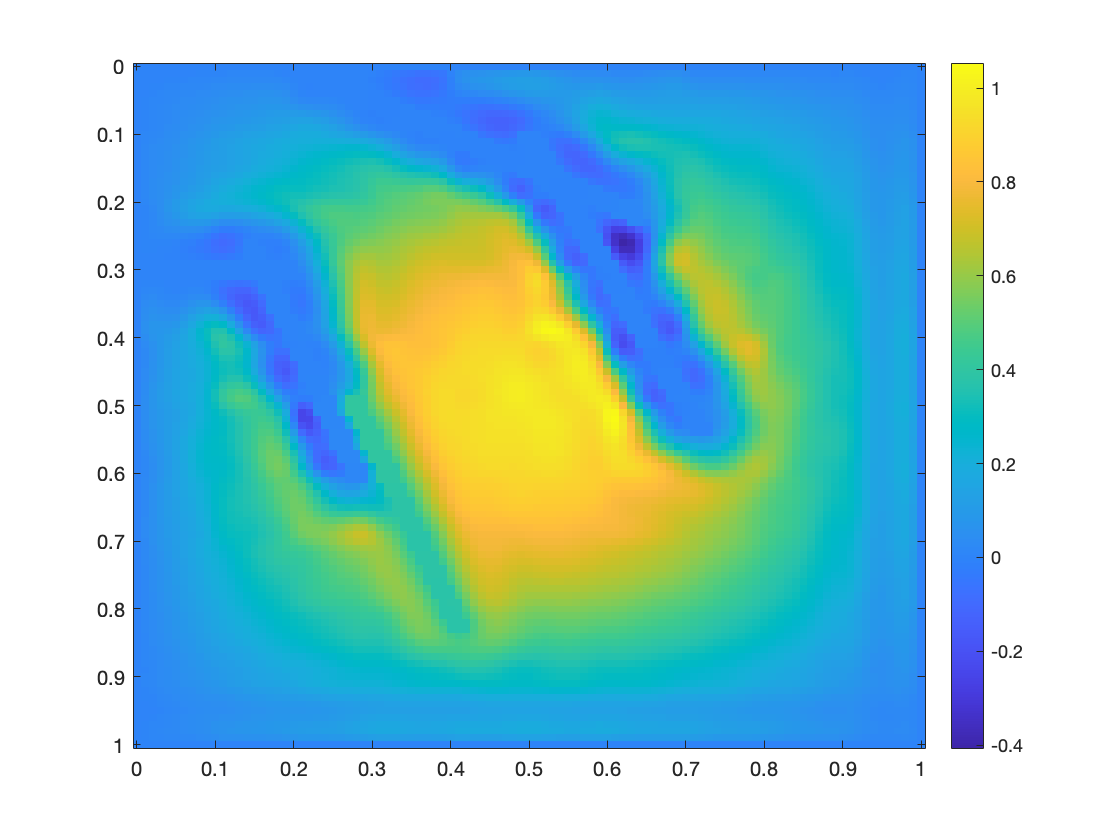}
\includegraphics[width = 5.2cm]{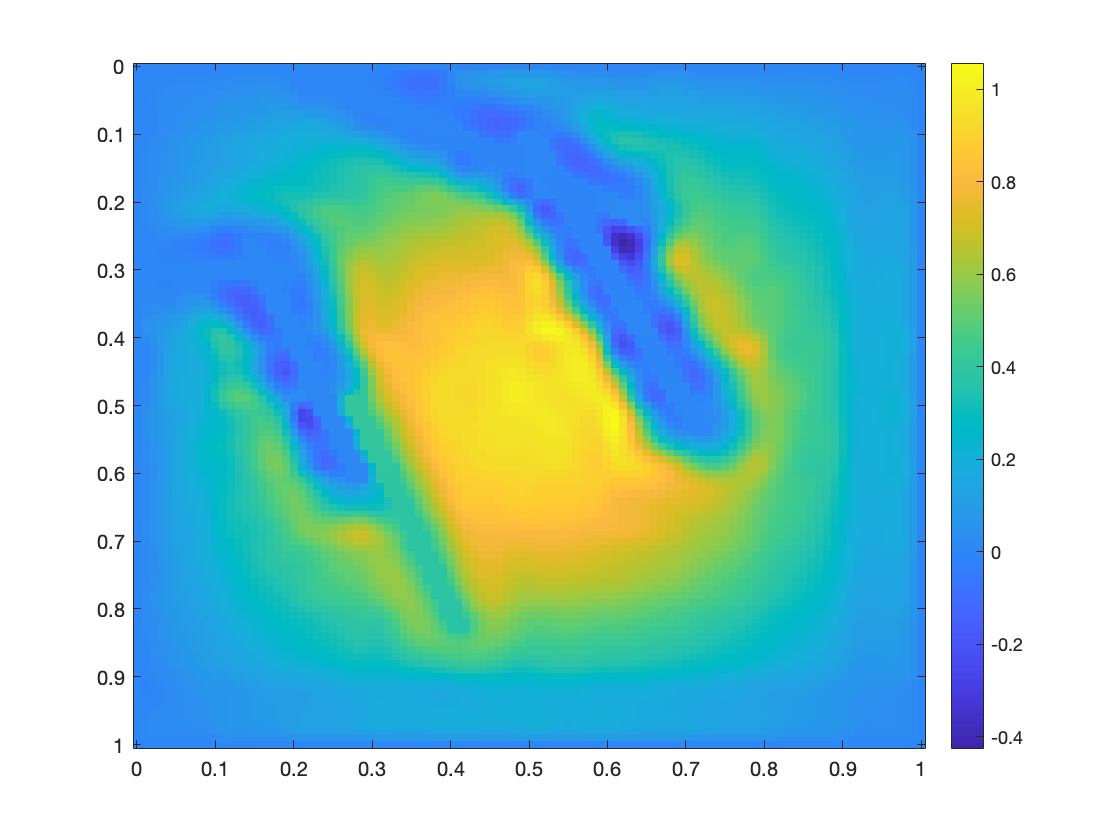}
\caption{Solution profiles at terminal time in Example \ref{exp:1}. Left: Reference solution. 
Middle: Implicit CEM solution with additional bases. 
Right: Partially explicit CEM solution.}
\label{NRfigic2}
\end{figure}

\begin{figure}[H]
\centering
\includegraphics[width = 8cm]{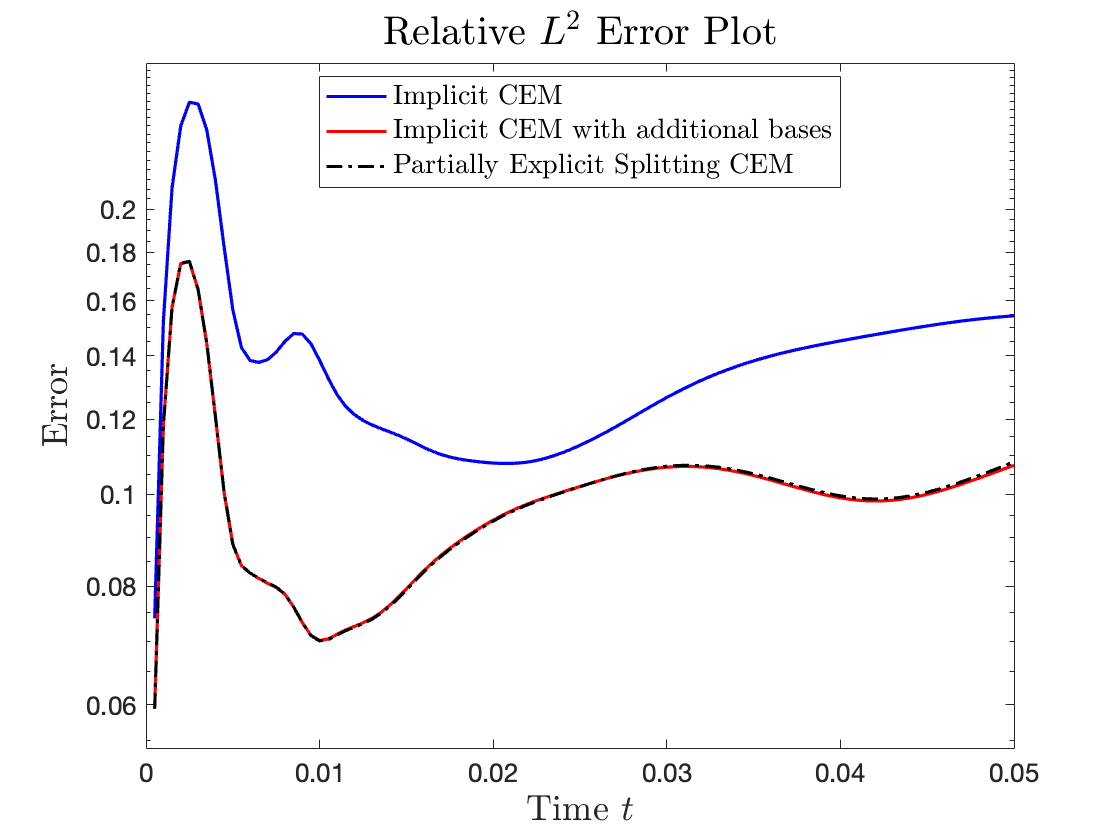}
\caption{Plot of relative $L^2$ error in Example \ref{exp:1}.}
\label{NRfigic3}
\end{figure}
\end{example}

\begin{figure}[H]
\centering
\includegraphics[width = 5.2cm]{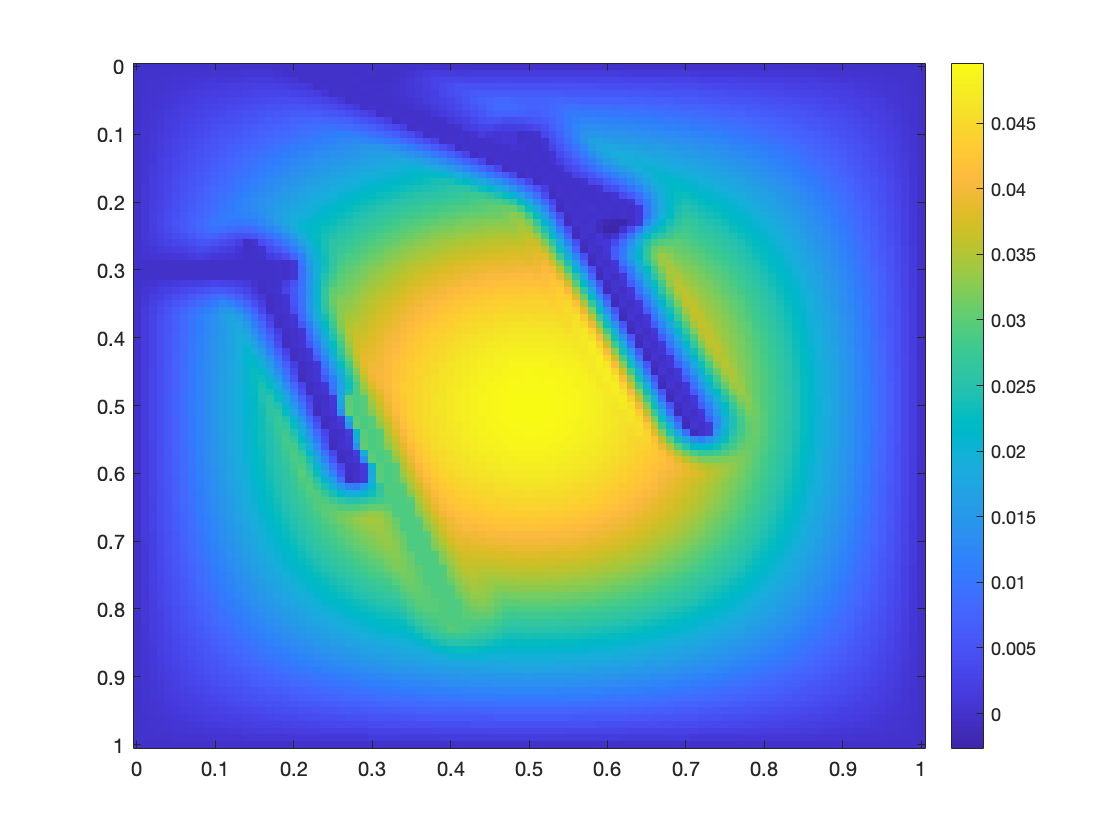}
\includegraphics[width = 5.2cm]{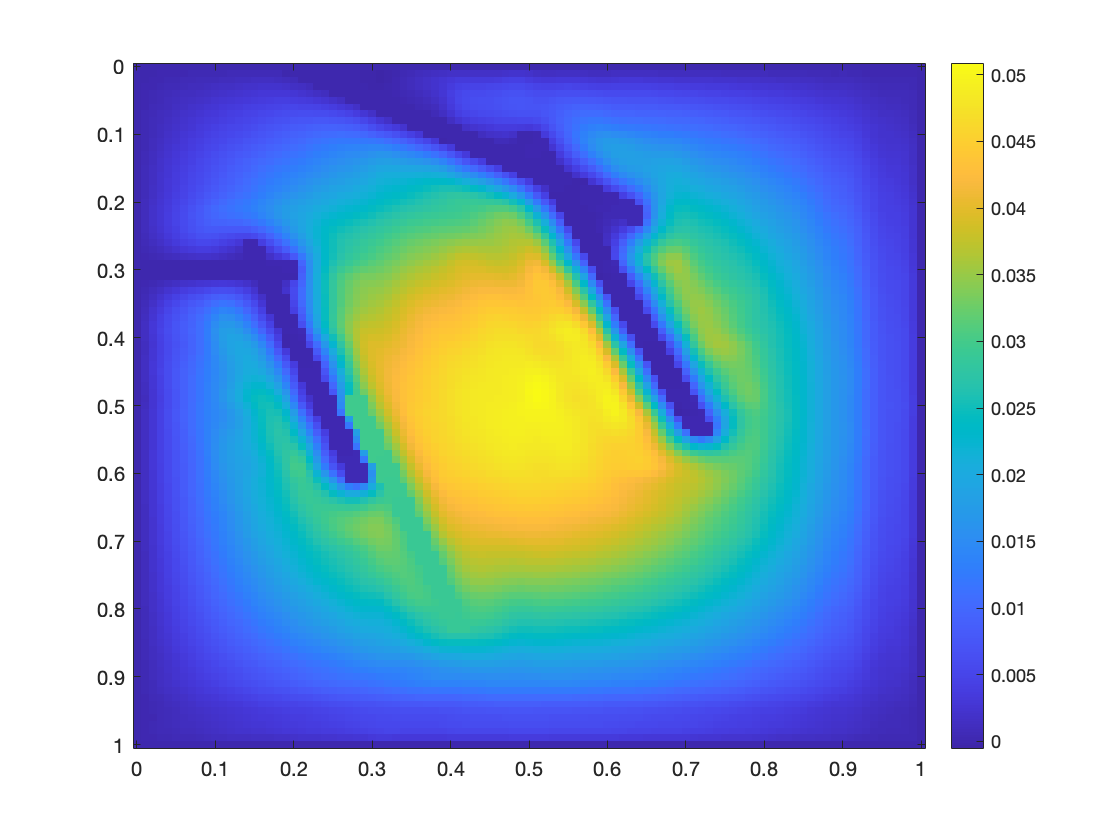}
\includegraphics[width = 5.2cm]{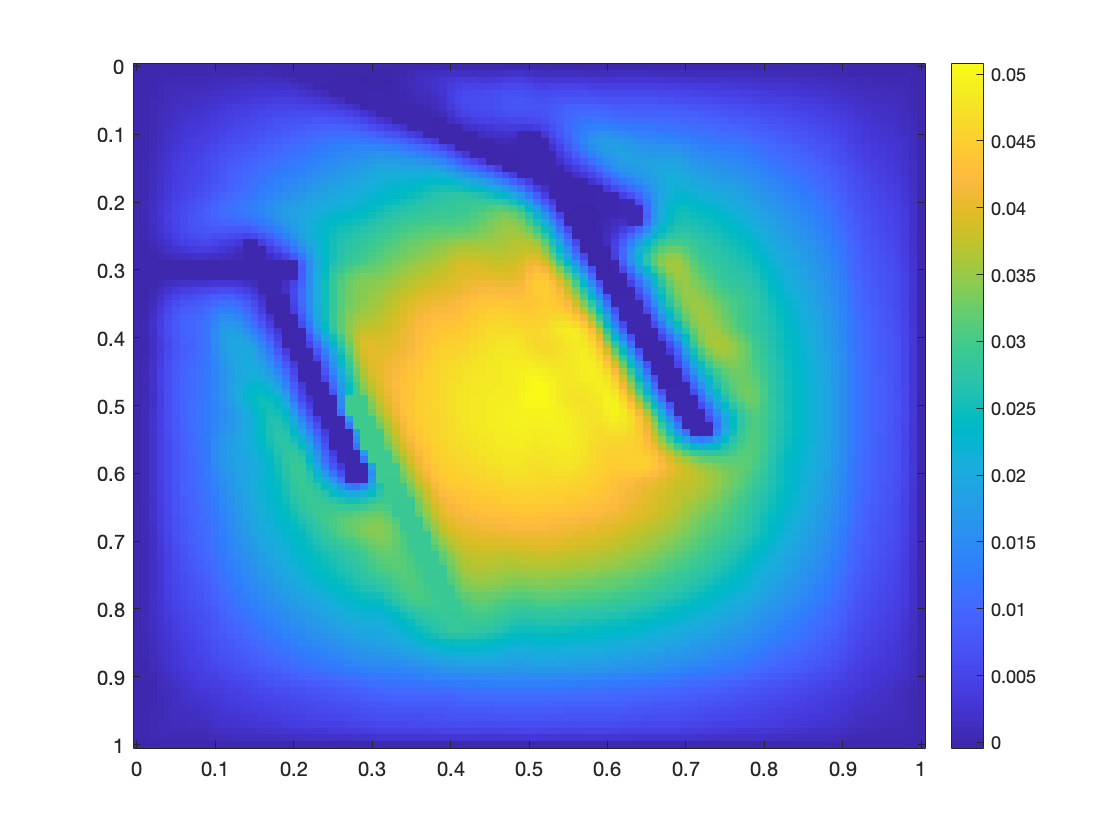}
\caption{Solution profiles at terminal time in Example \ref{exp:2}. Left: Reference solution. 
Middle: Implicit CEM solution with additional bases. 
Right: Partially explicit CEM solution.}
\label{NRfig2}
\end{figure}

\begin{figure}[H]
\centering
\includegraphics[width = 8cm]{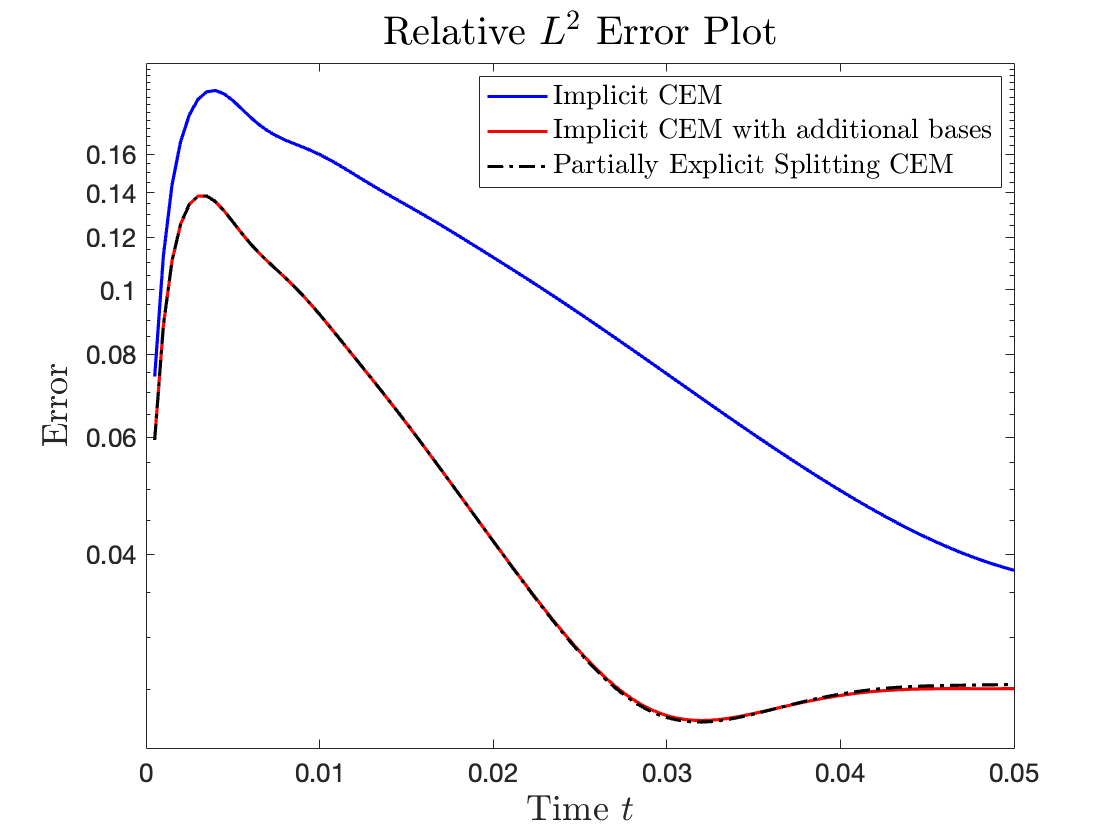}
\caption{Plot of relative $L^2$ error in Example \ref{exp:2}.}
\label{NRfig3}
\end{figure}

\begin{example} \label{exp:2}
In the second example, the permeability field $\kappa$ is the same as the one in Example \ref{exp:1}. We set the initial condition to be $u_0 = 0$. To avoid the solution $u$ being trivial, we 
 add a constant-in-time source term $g_0(x) = g_0(x_1,x_2) = \sin(\pi x_1) \sin(\pi x_2)$ for any $(x_1, x_2) \in \Omega$ to the right-hand side of the third and forth equations in \eqref{eqn:part-exp}. 
Let $T = 0.05$ and the temporal step size is $\Delta t = T/100 = 5 \times 10^{-4}$. 
%we show the permeability field $\kappa$ and the source term $g_0$ in Figure \ref{NRfig1}. We need to add the source term on the right hand side of the equations, otherwise the solution will be trivial. 

The profiles of the numerical solutions at the terminal time are sketched in Figure \ref{NRfig2}. 
The $L^2$ error curves against time are shown in Figure \ref{NRfig3}. 
%The reference solution is the implicit fine grid solution. 
%The blue, red and black curves in the plot represent the $L^2$ error for the implicit CEM solution, the implicit CEM solution with additional bases and the partially explicit solution, respectively. 
Similar to Example \ref{exp:1}, from the plot of error curves, we find that there is a considerable decrease in terms of $L^2$ error when we include $V_{H}^{2}$. Moreover, we notice that the curves for the implicit CEM solution with additional bases and the partially explicit solution nearly coincide, which implies that they have similar accuracy. In those settings, the $L^2$ error at the terminal time is about $1.9\%$. 
\end{example}

\begin{comment}
\begin{figure}[H]
\centering
\includegraphics[width = 7.5cm]{1a}
\includegraphics[width = 7.5cm]{1f}
\caption{Left: $\kappa$. Right: $g_0$.}
\label{NRfig1}
\end{figure}
\end{comment}

\begin{figure}[H]
\centering

\includegraphics[width =8cm]{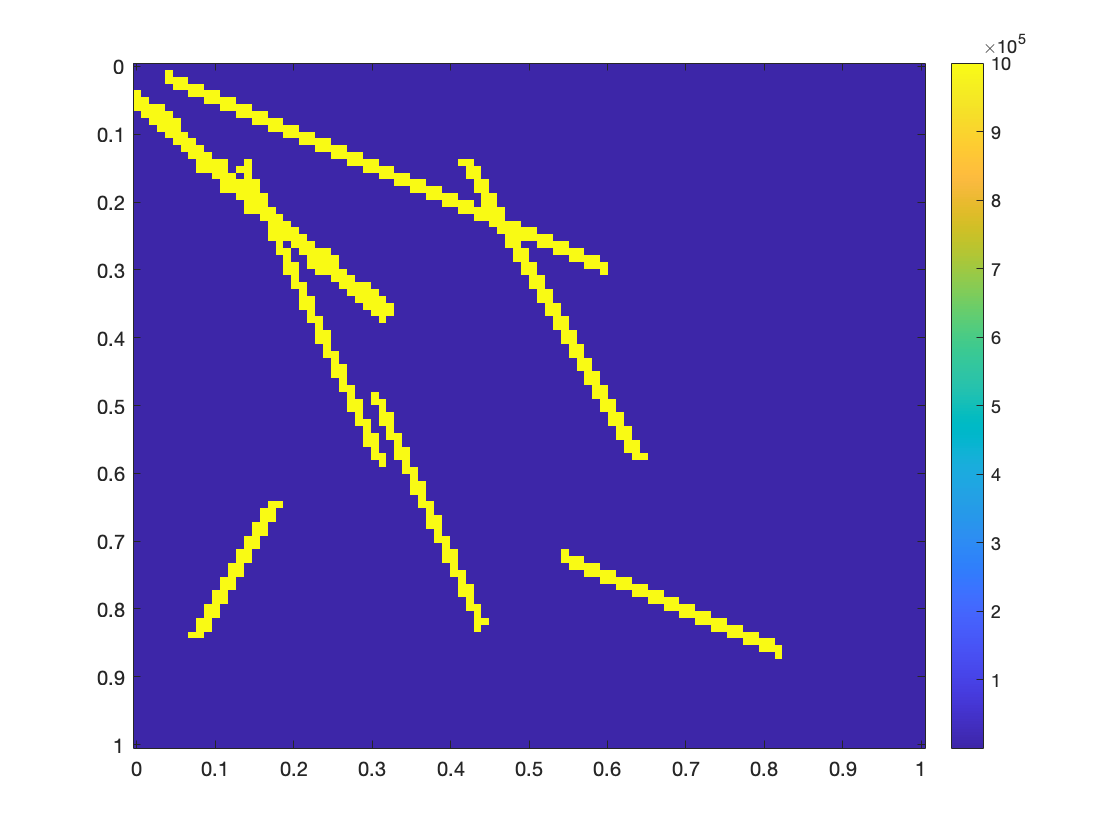}

\caption{Permeability in Example \ref{exp:3}.}
\label{NRfig4}
\end{figure}

\begin{figure}[H]
\centering
\includegraphics[width = 5.2cm]{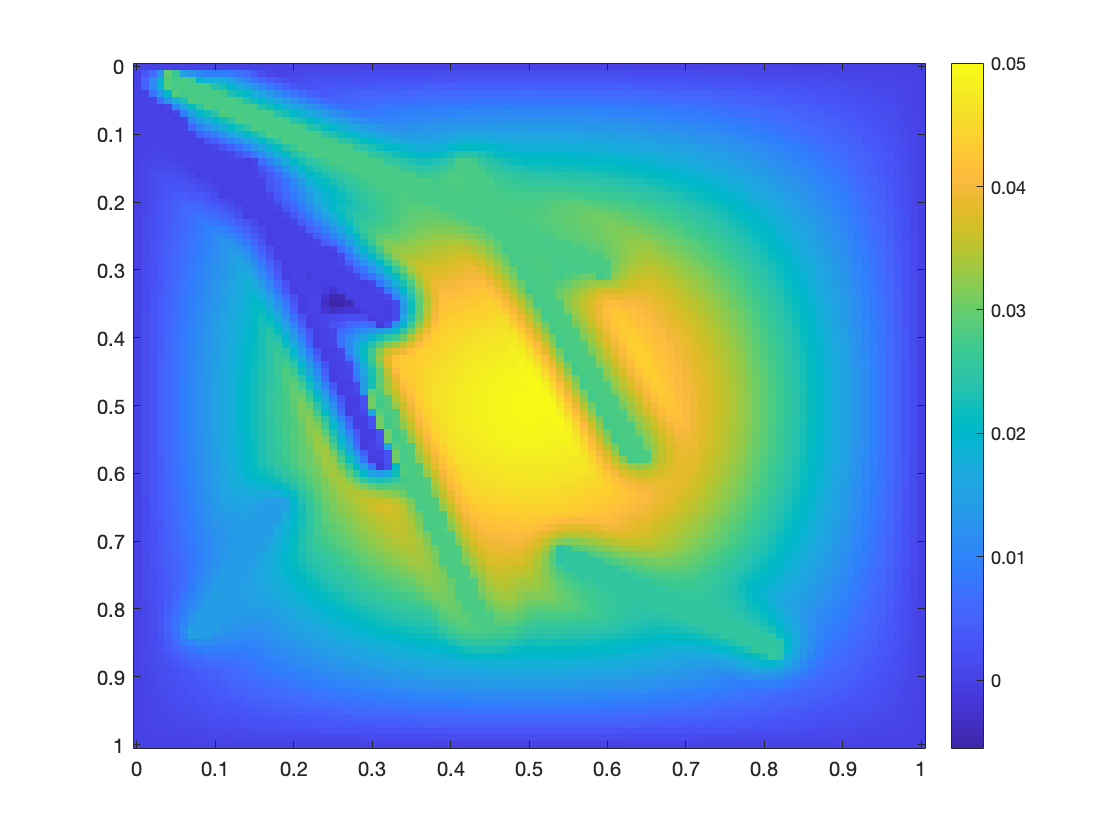}
\includegraphics[width = 5.2cm]{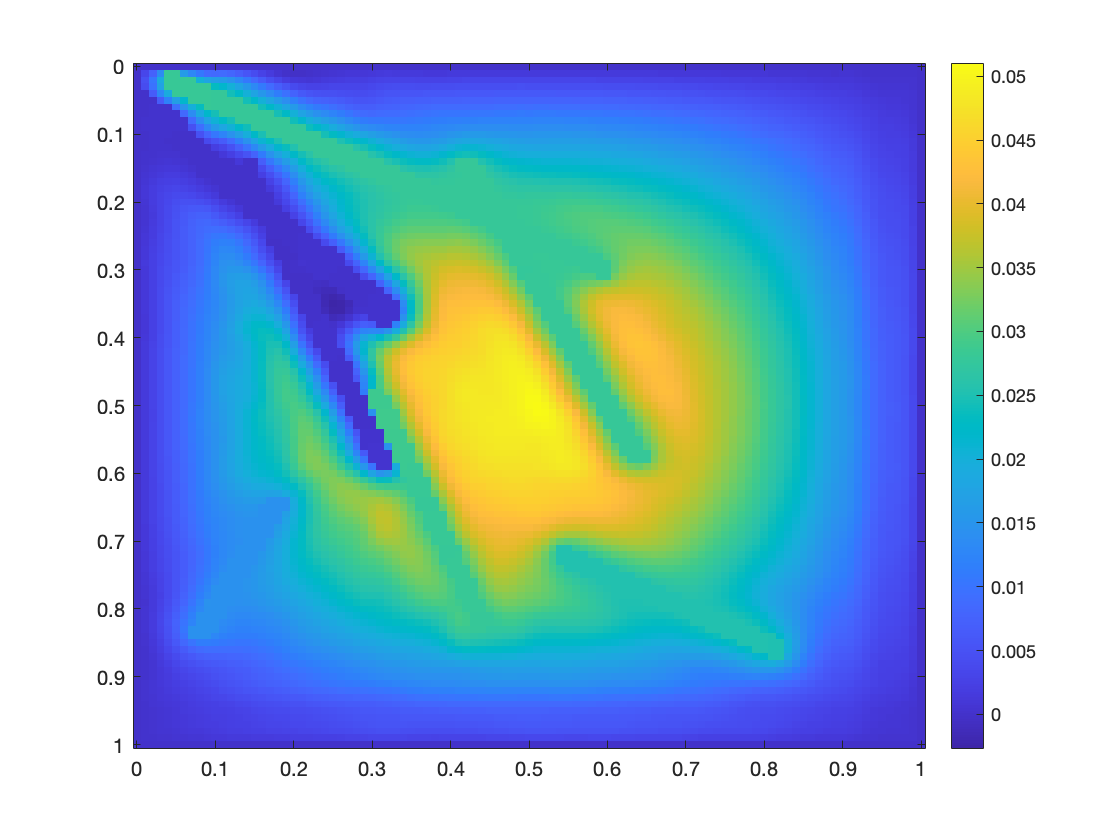}
\includegraphics[width = 5.2cm]{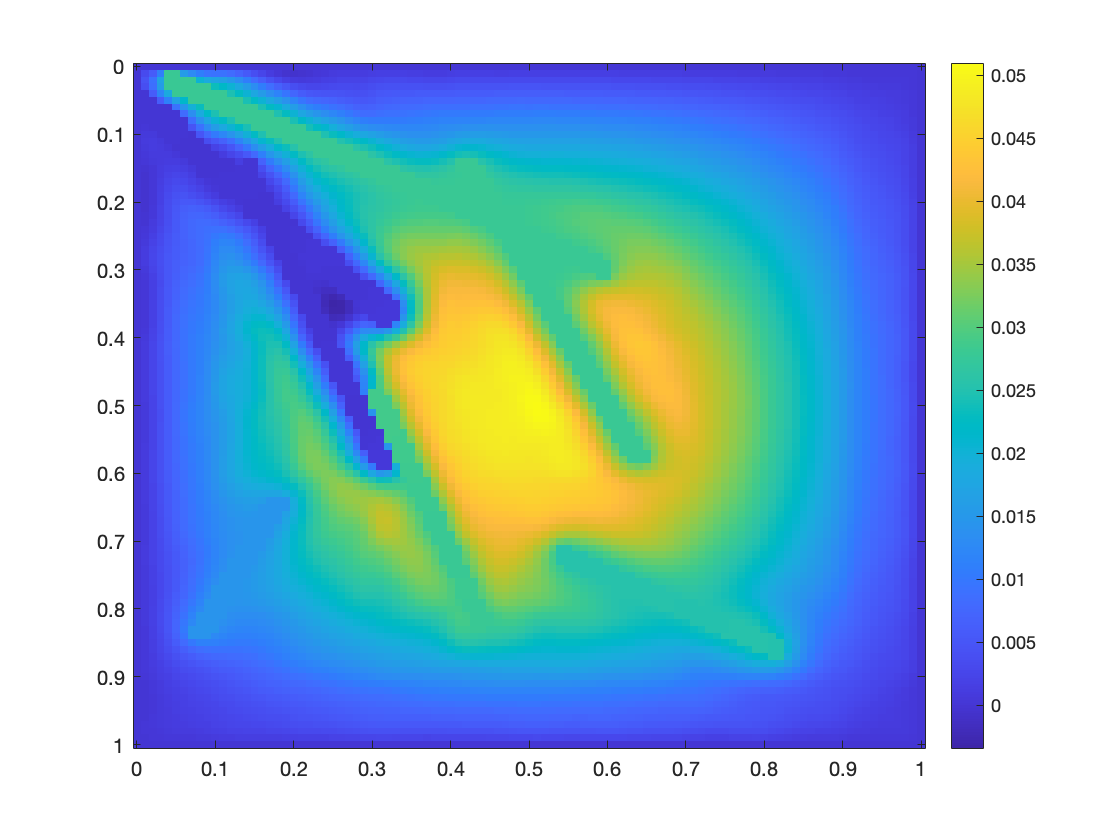}
\caption{Solution profiles at terminal time in Example \ref{exp:3}. Left: Reference solution. 
Middle: Implicit CEM solution with additional bases. 
Right: Partially explicit CEM solution.}
\label{NRfig5}
\end{figure}

\begin{figure}[H]
\centering
\includegraphics[width =8cm]{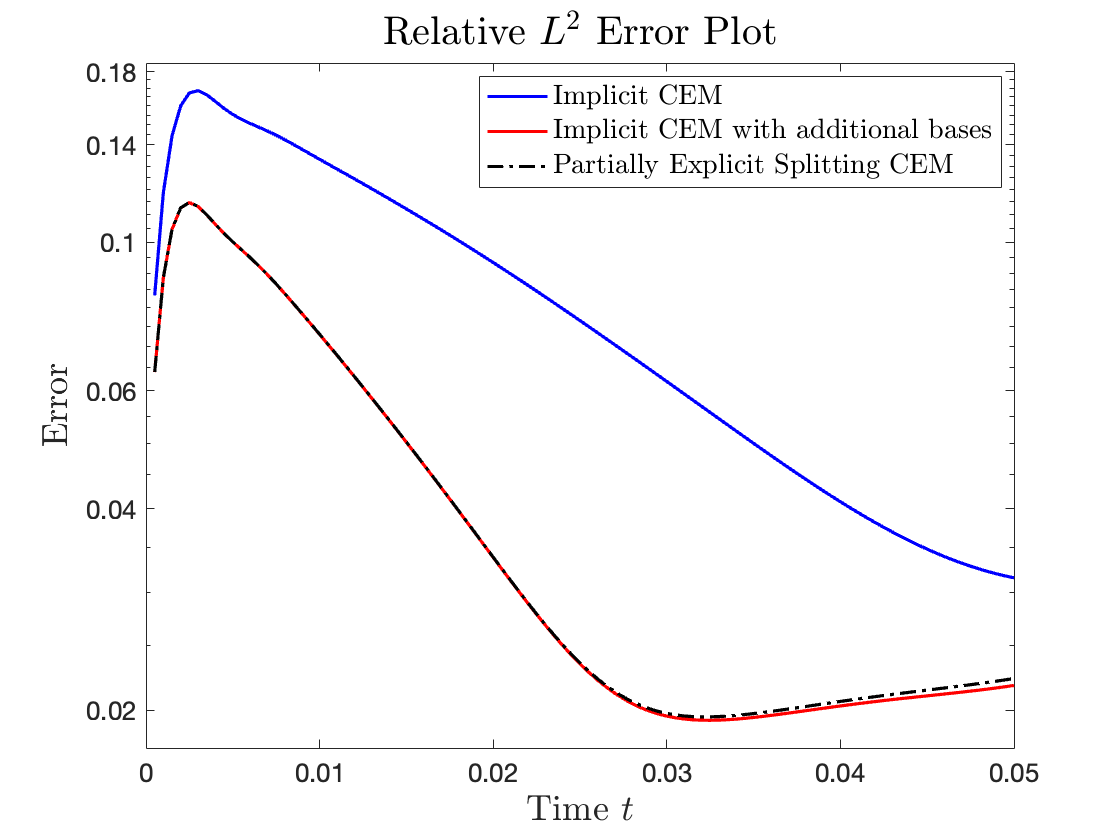}
\caption{Plot of relative $L^2$ error in Example \ref{exp:3}.}
\label{NRfig6}
\end{figure}

\begin{example}
\label{exp:3}
In the third example, we take the problem setting from Example \ref{exp:2} except the permeability. The permeability field for this case contains more {\it channels} and it is depicted in Figure \ref{NRfig4}. 
%The source term $g_0$ is the same as in Example \ref{exp:2}. 
%The implicit fine grid solution, the implicit CEM solution with additional bases and the partially explicit solution are presented in Figure \ref{NRfig5}. 
The solutions profiles are plotted in Figure \ref{NRfig5}; and the relative $L^2$ error plot is shown in Figure \ref{NRfig6}. 
%  We observe that the $L^2$ error for the implicit CEM solution with additional bases and the partially explicit solution are nearly the same. Therefore, the partially explicit scheme can obtain similar accuracy as the implicit CEM scheme with additional bases.   
In this case, the $L^2$ error at the terminal time is about $2.2\%$ using the partially explicit splitting scheme, which is comparable to the case using the implicit method with additional basis functions. 
This demonstrates the effectiveness and efficiency of the proposed partially explicit temporal discretization with additional basis functions from $V_H^2$. 
\end{example}

\section{Conclusion} \label{sec:conclusion}
In this work, we propose
 dememorization technique for a class of convection-diffusion equations
with memory effects.
These macroscopic  equations arise as a result of homogenization or upscaling
of transport equations (without memory terms) in heterogeneous media.
Because of transport at the microscales, the upscaled equations
contain memory terms.  
The dememorization technique introduced in the paper
 allows removing the term with memory effect and simplifying
the computations.
The dememorized equations differ from the original micro-scale 
equations. 
For the numerical discretization, we consider 
both implicit and partially explicit
splitting methods within the framework of CEM-GMsFEM. 
The latter scheme was previously introduced for 
problems in multiscale media with
high-contrast properties, which had been shown to be effective for such category of problems. 
Numerical results were presented that demonstrate the effectiveness and efficiency of the partially explicit schemes applying to the 
dememorized system of equations.

\bibliographystyle{abbrv}
\bibliography{references}

\appendix 
\section{Derivation of Macroscopic Equations}
\label{app:perturb}
%\marginpar{Change notations}
In this section, we derive the macro-scale equation with a nonlocal memory term from the micro-scale transport equation. 
We consider two cases. The first one is upscaling the transport equation with a perturbation perspective; while the second case deals with the flow transport in a special type of porous medium. 

\subsection{Perturbation Approach}
Consider the following transport equation in a heterogeneous medium 
\begin{eqnarray}
u_t^\veps + a^\veps \cdot \nabla u^\veps= -ku^\veps.
\label{eqn:pert}
\end{eqnarray}
In general, the velocity field $a^\veps : \mathbb{R}^d \to \mathbb{R}^d$ contains features at micro-scale driven by some hidden parameter $\veps$. 
Let $a^\veps=\vz +\vd$ and  $u^\veps=\cz +\cd$, where $\vd$ and $\cd$ are some {\it small} perturbations from the mean values.
Here, we can think of the case when $a^\veps$ (resp. $u^\veps$) is a {\it random} velocity field with $\vz$ (resp. $\cz$) being the {\it expectation} of the random field. 
Then, the equation \eqref{eqn:pert} becomes 
\begin{eqnarray}
  	(\cz +\cd)_t +({\vz +\vd}) \cdot \nabla(\cz +\cd) = -k(\cz +\cd).
	\label{eq:pertubation}
\end{eqnarray}
Taking the expectation on both sides of \eqref{eq:pertubation}, we have 
\begin{eqnarray}
  	\cz_t + \vz\cdot \nabla \cz + \overline{\vd\cdot\nabla \cd} = -k\cz . \label{eq:averaged}
\end{eqnarray}
Subtracting \eqref{eq:averaged} from \eqref{eq:pertubation}, we get
$$   	\cd_t +\vz\cdot \nabla \cd +\vd \cdot \nabla \cz + \underbrace{\vd \cdot \nabla  \cd - \overline{\vd \cdot \nabla \cd}}_{\text{high order terms}} = -k\cd $$
or equivalently, 
\begin{eqnarray}
  	\cd_t +\vz\cdot \nabla \cd + k\cd = - \vd \cdot \nabla \cz  + (\text{high order terms}).
	\label{eq:pertubsubsaver}
\end{eqnarray}
In general, the fluctuation $\cd$ is a function of $x$, $t$, and some hidden random variables; thus, we simply write $\cd = \cd (x,t)$. Notice that for any fixed $t>0$, we define $\varphi = \varphi(x,s;t)$ such that 
$$ \varphi (x,s;t) = \cd (\tilde x(x,s;t), s) \quad \text{for any} ~ s \in (0,t),$$
where $\tilde x = \tilde x(x,s;t)$ is a trajectory which satisfies 
$$ \frac{\partial \tilde x}{\partial s} = \vz, \quad \tilde x\vert_{s=t} = x.$$
%Next, we note that $$\frac{d}{ds}\cd (x(t,s),s) = \cd_t +\vz\cdot \nabla \cd,$$ where $x(t,s)$ is defined as above.
In the case of $\vz = \vz(x)$, we can rewrite the trajectory as $\tilde x = x - \vz (t-s)$. Then, the equation \eqref{eq:pertubsubsaver} becomes 
$$ \varphi_s  + k \varphi = - \vd \cdot \nabla \cz + (\text{high order terms}).$$
Assume $\varphi \vert_{s=0}= 0$. We integrate (with respect to $s$) the above equation over interval $(0,t)$, multiply the result by $\vd_i (x)$ (where $\vd = (\vd_i)_{i=1}^d$), and take average with respect to the randomness. This gives 
$$
	\overline{\vd_i (x) \cd } + k\int_0^t \overline{\vd_i \cd (\tilde x ,s) }\, ds 
	= -\overline{\vd_i (x) \int_0^t\vd_j (\tilde x){\partial \cz\over \partial x_j}(\tilde x, s)\, ds} 
	+ \text{(high order terms)}.
$$
Here, the Einstein summation convention is used for the index $j$. 
Neglecting the high order terms, this implies that $q = \overline{\vd_i \cd}$ satisfies the following ODE
$$ q_t + kq = -\overline{\vd_i \vd_j {\partial \cz\over \partial x_j}} .$$
One can solve for $q$ and it implies that %\textbf{Why $\mathbf{\tilde{a}_j(\tilde{x})}$?}
\begin{eqnarray}
\overline{\vd_i \cd} = - \overline{\vd_i \int_0^t e^{k(t-s)} \vd_j (x-\vz (t-s)) \frac{\partial \cz}{\partial x_j} (x-\vz (t-s))\, ds}.
\label{eq:averudcd}
\end{eqnarray}
Assume that $\vd$ is divergence-free. Substituting \eqref{eq:averudcd} into \eqref{eq:averaged}, we get 
\begin{eqnarray}
  	\cz_t + \vz\cdot \nabla \cz +k\cz  =
 {\partial \over \partial x_i}  \int_0^te^{-k(t-s)} \overline{\vd_i (x)  \vd_j (x- \vz(t-s))  {\partial \cz \over \partial x_j} (x-\vz(t-s))} \, ds. \label{eq:averaged2}
\end{eqnarray}
We remark that the Einstein summation convention is used for indices $i$ and $j$ in \eqref{eq:averaged2}. This gives the macroscopic transport equation for $\cz$ under the average velocity field $\vz$ with a memory term on the right-hand side.

%\appendix
\subsection{Upscaling for Layered Media}
\label{app:hom}
%To discuss the relation to homogenization, we consider two cases. 
In this case, we consider the transport equation
in a layered medium and seek the solution $u= u(x,t)$ such that 
$$
u_t + a \left ({x_2\over \varepsilon}\right)\cdot {\partial u\over \partial x_2}  =0,
$$
with the initial condition $u(x,0)=H(x_1)$ for any $x = (x_1, x_2)$ in a bounded domain, where $H(\cdot)$ is the single variable Heaviside function. The $a(\cdot)$ is a given scalar function describing the velocity along the $x_2$ direction. 
We are interested in the homogenized solution given by the average along the $x_2$ direction as follows: 
$$
\overline{u}=\int u(x,t) dx_2.
$$
We remark that the above integral should be understood in the sense of average along the direction of $x_2$. 
%\marginpar{Please, change $v$ to $a$ and $v_i$ to $a_i$ and $c$ to $u$.}
We present a discrete case assuming that $a$ takes
values $a_i$ in the $i$-th layer that has a width $m_i$ (with $n$ being total number of layers), i.e.,
$$
a\left (x_2\over \varepsilon\right) = a_i, \quad \text{if} ~ y_{i-1}\leq \frac{x_2}{\varepsilon} < y_i, \quad m_i = y_{i} - y_{i-1}
$$
for $i \in \{ 1, \cdots, n \}$. 
Then, the averaged solution 
can be written as  
$$\overline{u}(x,t)=\sum_{i=1}^n m_i H(x-a_i t),$$
and the homogenized solution is given by 
$${\partial \overline{u}\over \partial t} + \overline{a}{\partial \overline{u}\over \partial
x_1}=\sum_{i=1}^{n-1} \int_0^t \beta_i
{\partial^2 \overline{u}\over \partial x_1^2} (x_1-a_i(t-\tau),\tau)  d\tau,
$$
where $\vz$, $\beta_i$, and $u_i$ ($i=1,\dots,n-1$) satisfy
\begin{eqnarray*}
\begin{split}
\sum_{k=1}^{n}   {m_k\over u_i-a_k}&=0 & \quad i=1,\cdots,n-1, \\
\sum_{i=1}^{n-1}  {\beta_i\over u_i-a_k} &= (\overline{a}-a_k) & \quad  k=1,\cdots,n,  \\
\overline{a}& =\sum_{i=1}^n m_i a_i.&
\end{split}
\end{eqnarray*}
See \cite{efendiev2000exact} for details. 
Note that one can show that $\beta_i$'s and $u_i$'s exist and are unique. Moreover, they %depend only on one point statistics of the media, $M_l=\sum_{i=1}^n m_i v_i^l$, and 
have the following properties:
\begin{enumerate}
\item $a_1\leq u_1\leq a_2\leq\dots\leq u_{n-1}\leq a_n$; 

\item $\sum_{i=1}^{n-1}\beta_i=var(a)$, where $var(a)$ denotes
the variance of the velocity field and is given by
$var(a)=\sum_{i=1}^n m_ia_i^2-\left (\sum_{i=1}^n m_ia_i\right)^2$. 
\end{enumerate}

\begin{comment}
\subsection{Velocity with small perturbations}
%\marginpar{Change, $c$ to $u$, $v$ to $a$}
We assume that 
the velocity is a stochastic field, $a(x,\omega)$ with the mean $\overline{v}$
\[
a = \overline{a} + a',
\]
where $\overline{a}=\mathbb{E}_\omega(a)$ is the expectation with respect to some probability measure. In this case, we can also express 
\[
u = \overline{u} + u'.
\]
By substituting these perturbation arguments to the main equation \textbf{which equations?} and
taking the average, we arrive at
\[
\overline{u}_t + \overline{a} \cdot \nabla \overline{u}=
{\partial\over \partial x_i} \int_0^t \overline{a_j'(x)a_j'(x(t,s))}{\partial \over \partial x_j} \overline{u}(x(t,s),s) ds,
\]
where $x(t,s)$ ($s<t$) is trajectory such that $x(t,t)=x$ and
\[
{dx\over ds}  = \overline{a}.
\]
In the case of constant $\overline{a}=\lambda$,
\[
x(t,s) = x-\lambda(t-s).
\]
\end{comment}

\section{The Stability of Coupled System}
\label{app:stable}
We consider the case of $M=1$ in \eqref{eqn:model-coupled} for simplicity. The general case with $M>1$ can be derived similarly. The equations become
\begin{eqnarray}
\begin{split}
\partial_{t}v-\widetilde{a}\cdot\nabla v+\beta v & =u,\\
\partial_{t}u+a\cdot\nabla u & =\nabla\cdot\kappa\nabla v,
\end{split}
\label{eqn:csys}
\end{eqnarray}
with $u|_{\partial\Omega}=0$. 
One can also consider other boundary conditions rather than the homogeneous Dirichlet type. In this case, one has to assume that $\widetilde {a} \cdot \mathbf{n}_{\partial \Omega} = 0$. 
% and $\widetilde{a}\cdot n_{\partial\Omega}=0$\textbf{if we assume $w\in H^1_0$, then we don't need this condition.}.
For simplicity, we assume that the spatial dimension is $d = 2$. 

In the following, we consider the case when the permeability tensor $\kappa$ has the form 
$$\kappa=\kappa_{11}\widetilde{a}\otimes \widetilde{a}+\kappa_{22}\widetilde{a}^{\perp}\otimes \widetilde{a}^{\perp}+\kappa_{12}(\widetilde{a}^{\perp}\otimes \widetilde{a}+\widetilde{a}\otimes \widetilde{a}^{\perp}),$$
%where $\widetilde{a}^{\perp}\cdot \widetilde{a}=0$ and $\widetilde{a}^{\perp}\cdot \widetilde{a}^{\perp}=1$ with $\widetilde{a}(x)\neq 0$ for all $x$. 
where $\widetilde{a}$ and $\widetilde{a}^{\perp}$ are divergence free and $\widetilde{a} (x) \neq 0$ for almost all $x \in \Omega$. 
Here, $\kappa_{11}$, $\kappa_{22}$, and $\kappa_{12}$ are some heterogeneous scalar functions. 
We also write $\kappa_{21} = \kappa_{12}$. 
By applying $\nabla\cdot \kappa \nabla$ to the first equation, we have
\begin{eqnarray*}
\begin{split}
\partial_{t}\nabla\cdot(\kappa\nabla v)-\nabla\cdot(\kappa\nabla(\widetilde{a} \cdot \nabla v))+\beta \nabla\cdot(\kappa \nabla v) & =\nabla\cdot(\kappa \nabla u),\\
\partial_{t}u+a \cdot\nabla u & =\nabla\cdot\kappa\nabla v.
\end{split}
\end{eqnarray*}
Testing the first equation with $w\in H^1(\Omega)$ and using integration by parts, we obtain
$$
\int_{\Omega}\Big(\partial_{t}(\kappa\cdot\nabla v)-\kappa\nabla(\widetilde{a}\cdot\nabla v)+\beta (\kappa\nabla v)\Big)\cdot\nabla w ~ dx=\int_\Omega (\kappa\nabla u)\cdot\nabla w ~ dx \quad \forall w\in H^{1}(\Omega).
$$
Next, we define the following scalar functions such that 
$$
\tilde{v}_1 := \widetilde{a}\cdot\nabla v, \quad \tilde{v}_2 := \widetilde{a}^{\perp}\cdot\nabla v, \quad \tilde{w}_1 := \widetilde{a}\cdot\nabla w, \quad \text{and} \quad \tilde{w}_2 := \widetilde{a}^{\perp}\cdot\nabla w.
$$
Thus, we have 
\begin{eqnarray*}
\begin{split}
(\kappa\nabla v)\cdot\nabla w & =\sum_{1 \leq i,j\leq 2}\kappa_{ij}\tilde{v}_{i}\tilde{w}_{j},\\
(\kappa\nabla(\widetilde{a}\cdot\nabla v))\cdot\nabla w&=\sum_{i=1}^2 \tilde{w}_{i} \left [ \kappa_{1i}\tilde{a}+  \kappa_{2i}\tilde{a}^{\perp} \right ] \cdot\nabla\tilde{v}_{1}.
\end{split}
\end{eqnarray*}
We define a tensor $C = (C_{ijk})_{1 \leq i, j,k \leq 2}$ by 
$$
C_{ijk}:=\begin{cases}
\tilde{\kappa}_{1i}\widetilde{a}_{j}+\tilde{\kappa}_{2i}\widetilde{a}_{j}^{\perp} & \text{if} ~ k=1,\\
0 & \text{if} ~ k=2.
\end{cases}
$$
Then, we have 
$$
(\kappa\nabla(\widetilde{a}\cdot\nabla v))\cdot\nabla w = \sum_{i=1}^2 \tilde{w}_{i} \left [ \kappa_{1i}\widetilde{a}+  \kappa_{2i}\widetilde{a}^{\perp} \right ] \cdot\nabla\tilde{v}_{1} = \sum_{1 \leq i,j,k\leq 2} \tilde{w}_{i}C_{ijk}\partial_{j}\tilde{v}_{k}.
$$
As a result, we have
\begin{eqnarray*}
\begin{split}
\sum_{1 \leq i,j,k\leq 2} \int_{\Omega}\tilde{w}_{i}C_{ijk}\partial_{j}\tilde{v}_{k}  ~ dx&= \sum_{1 \leq i,j,k\leq 2} \int_{\Omega}\tilde{w}_{i}\left ( \cfrac{C_{ijk}+C_{kji}}{2} \right )\partial_{j}\tilde{v}_{k} ~ dx +\int_{\Omega}\tilde{w}_{i}\left (\cfrac{C_{ijk}-C_{kji}}{2}\right )\partial_{j}\tilde{v}_{k}~dx \\
& =  \sum_{1 \leq i,j,k\leq 2} \int_{\Omega}\tilde{w}_{i}\left ( \cfrac{C_{ijk}+C_{kji}}{2} \right )\partial_{j}\tilde{v}_{k} ~ dx
\end{split}
\end{eqnarray*}
since 
$$ \sum_{1 \leq i,j,k \leq 2} C_{ijk} - C_{kji} = 0.$$
Moreover, we have 
\begin{align*}
&\int_{\Omega}\tilde{v}_{i}\left (\cfrac{C_{ijk}+C_{kji}}{2}\right)\partial_{j}\tilde{v}_{k} ~ dx\\
 =&\cfrac{1}{2}\int_{\Omega}\tilde{v}_{i}\cfrac{C_{ijk}+C_{kji}}{2}\partial_{j}\tilde{v}_{k} ~ dx -\cfrac{1}{2}\left [\int_{\Omega}(\partial_{j}\tilde{v}_{i})\cfrac{C_{ijk}+C_{kji}}{2}\tilde{v}_{k} ~ dx + \int_\Omega \tilde{v}_{i}\partial_{j} \left ( \cfrac{C_{ijk}+C_{kji}}{2} \right )\tilde{v}_{k} ~ dx \right ]\\
 =&-\cfrac{1}{2}\int_{\Omega}\tilde{v}_{i} \partial_{j}\left ( \cfrac{C_{ijk}+C_{kji}}{2} \right )\tilde{v}_{k} ~ dx.
\end{align*}
Therefore, the system \eqref{eqn:csys} is stable if
$$ \beta \sum_{1 \leq i,j \leq 2} \kappa_{ij} + \frac{1}{4} \sum_{1 \leq i , j , k \leq 2}  (C_{ijk} + C_{kji}) \geq 0.$$
%$$\cfrac{1}{\beta}\partial_{j}\cfrac{C_{ijk}+C_{kji}}{2}<\tilde{\kappa}.$$

In particular, if $\kappa=\kappa_{11}\widetilde{a}\otimes \widetilde{a}$,  %and $\nabla\cdot a=\nabla\cdot\tilde{a}=0$,
%(\textbf{is $\nabla\cdot\tilde{a}=0$ needed?}) yes, it is required. 
we have %$\cfrac{C_{ijk}-C_{kji}}{2}=0$ and thus we have 
$$
\cfrac{1}{2}\Big(\|u(T)\|^{2}-\|u(0)\|^{2}\Big)=\int_{0}^{T}\int_{\Omega}u \left (  \partial_{t}u + a\cdot\nabla u \right ) dx dt=\int_{0}^{T}\int_{\Omega}u(\nabla\cdot\kappa\nabla v) ~ dx dt
$$
and 
$$
\int_{\Omega}\left ( v(\partial_{t}\nabla\cdot\kappa\nabla v) - (\tilde{a}\cdot\nabla v)\nabla\cdot\kappa\nabla v+ \beta v\nabla\cdot\kappa\nabla v \right ) dx =\int_{\Omega}v\nabla\cdot\kappa\nabla u ~ dx .
$$
Thus, integrating over $(0,T]$, we have 
$$\cfrac{1}{2}\left(\|\kappa_{11}^{\frac{1}{2}}(\tilde{a}\cdot\nabla v)(T)\|^{2}-\|\kappa_{11}^{\frac{1}{2}}(\tilde{a}\cdot\nabla v)(0)\|^{2}\right)+ \int_{0}^{T} \int_\Omega (\beta \kappa_{11} +\widetilde{a} \cdot \nabla\kappa_{11})(\widetilde{a}\cdot\nabla v)^{2}~ dxdt =  -\int_{0}^{T}\int_{\Omega}u(\nabla\cdot\kappa\nabla v) ~ dx dt.$$
%\begin{align*}
% & \cfrac{1}{2}\Big(\|\kappa_{11}^{\frac{1}{2}}(\tilde{a}\cdot\nabla v)(T)\|^{2}-\|\kappa_{11}^{\frac{1}{2}}(\tilde{a}\cdot\nabla v)(0)\|^{2}\Big)+ \int_{0}^{T} \int_\Omega (\beta \kappa_{11} +\nabla\kappa_{11}\cdot \widetilde{a})(\widetilde{a}\cdot\nabla v)^{2}~ dxdt\\
%= & -\int_{0}^{T}\int_{\Omega}u\nabla\cdot\kappa\nabla v
%\end{align*}
Therefore, we have 
$$\cfrac{1}{2}\Big(\|u(T)\|^{2}+\|\kappa_{11}^{\frac{1}{2}}(\tilde{a}\cdot\nabla v)(T)\|^{2}\Big)\leq\cfrac{1}{2}\Big(\|u(0)\|^{2}+\|\kappa_{11}^{\frac{1}{2}}(\tilde{a}\cdot\nabla v)(0)\|^{2}\Big)$$
%or 
%\[\cfrac{1}{2}\Big(\|u(T)\|^{2}+\|v(T)\|_{\kappa}^{2}\Big)\leq\cfrac{1}{2}\Big(\|u(0)\|^{2}+\|v(0)\|_{\kappa}^{2}\Big)\]
if $\beta \kappa_{11} + (\widetilde{a} \cdot \nabla\kappa_{11})\geq0$.
%We remark that the stable of $v$ can be given as 
%\[\cfrac{1}{2}\|v(T)\|^{2}+\beta \int_{0}^{T}\|v(t)\|^{2}=\int_{0}^{T}\int_{\Omega}uv+\cfrac{1}{2}\|v(0)\|^{2}\]
In this case, the stability of \eqref{eqn:csys} depends only on $\beta$, $\kappa_{11}$, and $\widetilde{a}$. 
%\begin{align*}
%\partial_{t}(\tilde{a}\cdot\nabla v)-\tilde{a}\cdot\nabla(\tilde{a}\cdot\nabla v)+\beta \tilde{a}\cdot\nabla v & =\tilde{a}\cdot\nabla u\\
%\partial_{t}u+a\cdot\nabla u & =\nabla\cdot\Big(\kappa_{11}\tilde{a}(\tilde{a}\cdot\nabla v)\Big)
%\end{align*}
%which is only depend on $u$ and $\tilde{a}\cdot\nabla v$.
%We remark that if $(\nabla\kappa_{11}\cdot \tilde{a})+\beta \kappa_{11}\leq0$,
%the equation should be unstable.

\section{Construction of the Ansatz Space}
\label{SpaceConstructSection}
%\textbf{need rephrasing}\\
In this section, we present the construction of the ansatz space $V_H = V_H^1 \oplus V_H^2$ that will be used for the spatial discretization. This ansatz space is based on the framework of the recently developed CEM-GMsFEM. 
For the ansatz space $W_H$, one can define it as the direct sum of $V_H$ and the span of the degrees of freedom corresponding to the in-flow boundary $\Gamma$. 
In the following, we define $V(S) := H_0^1(S)$ for a (nonempty) proper subset $S\subset \Omega$. In the following, we denote $V = H_0^1(\Omega)$. 

\subsection{The Implicit Ansatz Space}
\label{sec:cem}
In this section, we present the construction of the implicit ansatz space $V_H^1$. 
%The CEM method follows the framework of finite
%element methods. 
The construction of this space starts by solving a class of constrained energy minimization problems. 
Let $\mathcal{T}_{H}$
be a coarse grid partition of $\Omega$. 
Denote $N_e$ the total number of coarse elements.
For $K_{i}\in\mathcal{T}_{H}$,
we first have to build a collection of auxiliary bases in $V(K_{i})$.
Let $\{\chi_i\}_{i=1}^{N_c}$ be a set of partition of unity functions corresponding to an overlapping partition of the domain.
In each coarse element $K_i \in \mathcal{T}_H$, we solve the following eigenvalue problem:
\begin{eqnarray*}
\int_{K_i} \kappa \nabla \psi_j^{(i)} \cdot \nabla v = \lambda_j^{(i)} s_i ( \psi_j^{(i)},v) \quad
\forall v \in V(K_i),
\end{eqnarray*}
where 
$$s_i(u,v) = \int_{K_i} \tilde{\kappa} u v, \quad
\tilde{\kappa} := \kappa H^{-2} \; \text{or} \; 
\tilde{\kappa} := \kappa \sum_{i}\left|\nabla \chi_{i}\right|^{2}.$$
We rearrange and gather the $L_i$ eigenfunctions corresponding to the first $L_i$ smallest eigenvalues. 
Define the auxiliary space 
$$ V_{aux} := \bigoplus_{i=1}^{N_e}V_{aux}^{(i)} ,\quad  V_{aux}^{(i)}:=\text{span}\{\psi_{j}^{(i)}:1\leq j\leq L_{i}\}$$
and the projection operator $\Pi:L^{2}(\Omega)\to V_{aux}$ such that 
$$s(\Pi u,v)=s(u,v)\quad\forall v\in V_{aux}, \quad \text{where} ~ s(u,v):=\sum_{i=1}^{N_{e}}s_{i}(u|_{K_{i}},v|_{K_{i}}).$$ 
For an oversampling parameter $m \in \mathbb{N}$, we define $K_{i,m}$ to be an oversampling domain of $K_{i}$ as follows 
$$ K_{i,0} := K_i, \quad K_{i,m} := \bigcup \{ K \in \mathcal{T}_H : K \cap K_{i,m-1} \neq \emptyset \} \quad \text{for} ~ m \geq 1.$$
We simply denote $K_i^+ = K_{i,m}$ for some given oversampling parameter $m$. 
For each auxiliary basis $\psi_{j}^{(i)}$, we search for
a local basis function $\phi_{j}^{(i)}\in V(K_{i}^{+})$
such that 
\begin{eqnarray*}
\begin{split}
a(\phi_{j}^{(i)},v)+s(\mu_{j}^{(i)},v) & =0\quad &\forall v\in V(K_{i}^{+}),\\
s(\phi_{j}^{(i)},\nu) & =s(\psi_{j}^{(i)},\nu)\quad&\forall\nu\in V_{aux}(K_{i}^{+}),
\end{split}
\end{eqnarray*}
for some $\mu_{j}^{(i)} \in V_{aux}$, where $V_{aux}(K_i^+) := \displaystyle{\bigoplus_{K_j \subset K_i^+} V_{aux}^{(j)}}$. 
The implicit ansatz space $V_H^1$ is defined to be 
$$V_{H}^1 :=\text{span}\{\phi_{j}^{(i)}:\;1\leq i\leq N_{e},1\leq j\leq L_{i}\}.$$
%The CEM solution $u_{cem}$ is given by
%\begin{align*}
%(\cfrac{ \partial u_{cem} }{\partial t} ,v ) + ( f(u_{cem}),v) + (g(u_{cem}) , v ) = 0, \quad \forall v\in V_{cem}.
%\end{align*}
Let $\tilde{V}:= \{v\in V:  \Pi(v) = 0 \}$. Based on the construction of $V_H^1$, we have the property that $V = V_H^1 \perp_a \tilde V$. 

\subsection{The Explicit Ansatz Space}
In this section, we construct the explicit ansatz space $V_H^2 \subset \tilde V$. For each coarse element $K_{i}$, we consider the following class of eigenvalue problems: %We find eigenpairs $(\xi_{j}^{(i)},\gamma_{j}^{(i)})\in(V(K_{i})\cap\tilde{V})\times\mathbb{R}$ by solving
%For each $K_i$, after solving
find $\xi_{j}^{(i)}\in V(K_{i})\cap\tilde{V}$ and $\gamma_{j}^{(i)} \in \mathbb{R}$ such that 
$$
\int_{K_{i}}\kappa\nabla\xi_{j}^{(i)}\cdot\nabla v ~dx =\gamma_{j}^{(i)}\int_{K_{i}}\xi_{j}^{(i)}v ~dx  \quad  \forall v\in V(K_{i})\cap\tilde{V}.$$
%and rearranging the eigenvalues by $\gamma_1^{(i)}\leq \gamma_2^{(i)}\leq \cdots $., we rearrange and select the first $J_i$ eigenfunctions corresponding to the smallest $J_i$ eigenvalues. 
We define $V_{aux,2} := \text{span}\{\xi_j^{(i)} : 1\leq i \leq N_e, 1\leq j\leq J_i \}$.
For each $\xi_j^{(i)} \in V_{aux,2}$, we define $\zeta_{j}^{(i)} \in V(K_i^+)$ such
that for some $\mu_{j}^{(i),1} \in V_{aux}$, $ \mu_{j}^{(i),2} \in V_{aux,2}$, we have 
\begin{eqnarray*}
\begin{split}
a(\zeta_{j}^{(i)},v)+s(\mu_{j}^{(i),1},v)+ ( \mu_{j}^{(i),2},v) & =0 & \quad\forall v\in V(K_i^+), \\
s(\zeta_{j}^{(i)},\nu) & =0 &\quad\forall\nu\in V_{aux},\\
(\zeta_{j}^{(i)},\nu) & =( \xi_{j}^{(i)},\nu) & \quad\forall\nu\in V_{aux,2}.
\end{split}
\end{eqnarray*}
We define $V_{H}^2 :=\text{span}\{\zeta_{j}^{(i)}: 1\leq i \leq N_e, \;  1 \leq j\leq J_i\}$. Based on the construction, we have $\zeta_j^{(i)} \in \tilde V$ and thus $V_H^2 \subset \tilde V$. 

\section{Proof of Theorem \ref{thm:pe}}
\label{app:proof}
In the following, we omit the subscript $H$ and simply write $u = u_H$, $v = v_H$,  $u_i^n = u_{H,i}^n$, and $v_i^n = v_{H,i}^n$ to simplify the notations. 
Note that the first two equations in \eqref{eqn:part-exp} can be written as 
\begin{eqnarray*}
\begin{split}
\frac{1}{\Delta t} \left(v_{1}^{n+1} - v_1^n, \phi_1 \right) + \beta  (v_{1}^{n},\phi_1) & = (u_{1}^{n+1},\phi_1) & \quad \forall \phi_1 \in V_{H}^{1}, \\
\frac{1}{\Delta t} \left(v_{2}^{n+1} - v_2^n, \phi_2 \right) + \beta (v_{2}^{n},\phi_2) & =  (u_{2}^{n},\phi_2)  & \quad \forall \phi_2\in V_{H}^{2}.
\end{split}
\end{eqnarray*}

Taking $\psi_1 = u_{1}^{n+1}\in V_H^{1}$ and $\psi_2 = u_{2}^{n+1}\in V_H^{2}$ in \eqref{eqn:part-exp}, we obtain 
\begin{eqnarray*}
\begin{split}
\left (\cfrac{u_{1}^{n+1}-u_{1}^{n}}{\Delta t}+\cfrac{u_{2}^{n+1}-u_{2}^{n}}{\Delta t},u_{1}^{n+1}\right )+\mca(v_{1}^{n+1}+v_{2}^{n+1},u_{1}^{n+1}) & =0,\\
\left (\cfrac{u_{1}^{n+1}-u_{1}^{n}}{\Delta t}+\cfrac{u_{2}^{n+1}-u_{2}^{n}}{\Delta t},u_{2}^{n+1}\right )+\mca(v_{1}^{n+1}+v_{2}^{n+1},u_{2}^{n+1}) & =0,
\end{split}
\end{eqnarray*}
since $v^{n+1}=v_{1}^{n+1}+v_{2}^{n+1}$. Due to the construction of $V_H^1$ and $V_H^2$, we have $\mca(v_j^{n+1}, u_k^{n+1}) = 0$ for any $j \neq k$. Note that, taking $\phi_1 = v_1^{n+1}$, $\phi_2 = v_2^{n+1}$, and making use of the operator $\mcr$, we obtain 
\begin{eqnarray*}
\begin{split}
\mca(u_{1}^{n+1},v_1^{n+1}) & = \frac{1}{\Delta t} \mca(v_{1}^{n+1} - v_1^n ,v_1^{n+1}) + \beta \mca(v_{1}^{n+1},v_1^{n}),\\
\mca(u_{2}^{n},v_2^{n+1}) & = \frac{1}{\Delta t} \mca(v_{2}^{n+1} - v_2^n ,v_2^{n+1}) +\beta \mca(v_{2}^{n},v_2^{n+1}).
\end{split}
\end{eqnarray*}
Then, we have 
\begin{eqnarray*}
\begin{split}
& \quad \mca(v_{1}^{n+1}+v_{2}^{n+1},u_{1}^{n+1})+\mca(v_{1}^{n+1}+v_{2}^{n+1},u_{2}^{n+1}) \\
& = \mca(v_1^{n+1}, u_1^{n+1}) + \mca(v_2^{n+1}, u_2^{n}) + \mca(v_2^{n+1}, u_2^{n+1} - u_2^n) \\
%+ \mca(v_2^{n+1}, u_1^{n+1}) + \mca(v_1^{n+1}, u_2^{n+1}) \\
& =  \sum_{i=1}^2 \left [ \left ( \frac{1}{\Delta t} - \beta \right ) \mca(v_i^{n+1} , v_i^{n+1} - v_i^n) + \beta \mca(v_i^{n+1},v_i^{n+1}) \right ]  + \mca( v_2^{n+1}, u_2^{n+1} - u_2^n ). %+ \mca(v_2^{n+1}, u_1^{n+1}) + \mca(v_1^{n+1}, u_2^{n+1}).
\end{split}
\end{eqnarray*}
\begin{comment}
\begin{align*}
 & a(v_{1}^{n+1}+v_{2}^{n+1},u_{1}^{n+1})+a(v_{1}^{n+1}+v_{2}^{n+1},u_{2}^{n+1}) = a(v^{n+1},u_{1}^{n+1}+u_{2}^{n})+a(v^{n+1},u_{2}^{n+1}-u_{2}^{n})\\
= & -e^{-b\Delta t}\tilde{b}^{-1}a(v^{n},v^{n+1})+\tilde{b}^{-1}a(v^{n+1},v^{n+1})+a(v^{n+1},u_{2}^{n+1}-u_{2}^{n}).\\
= & e^{-b\Delta t}\tilde{b}^{-1}a(v^{n+1}-v^{n},v^{n+1})+(1-e^{-b\Delta t})\tilde{b}^{-1}a(v^{n+1},v^{n+1})+a(v^{n+1},u_{2}^{n+1}-u_{2}^{n}).
\end{align*}
\end{comment}
On the other hand, we have 
$$
\mca(v_i^{n+1} , v_i^{n+1} - v_i^n)=\cfrac{1}{2}\left (\|v_i^{n+1}\|_{\mca}^{2}-\|v_i^{n}\|_{\mca}^{2}+\|v_i^{n+1}-v_i^{n}\|_{\mca}^{2}\right) \quad \text{for any} ~ i \in \{ 1, 2\}
$$ 
%$$(1-e^{-b\Delta t})\tilde{b}^{-1}a(v^{n+1},v^{n+1})=b(\|v^{n+1}\|_{a}^{2}) \geq \cfrac{b}{2} (\|v^{n+1}\|_{a}^{2}),$$
%$$ \abs{\mca (v_j^{n+1}, u_k^{n+1})}  \leq \frac{1}{2} \left ( \norm{v_j^{n+1}}_a^2 +  \norm{u_k^{n+1}}_a^2 \right )\quad \text{for any} ~ j \neq k,$$
and 
$$|\mca(v_2^{n+1},u_{2}^{n+1}-u_{2}^{n})| \leq\|v_2^{n+1}\|_{\mca} \cdot \|u_{2}^{n+1}-u_{2}^{n}\|_{\mca} \leq\cfrac{\beta}{2}\|v_2^{n+1}\|_{\mca}^{2}+\cfrac{1}{2\beta}\|u_{2}^{n+1}-u_{2}^{n}\|_{\mca}^{2}.$$
Therefore, we have 
\begin{eqnarray}
\begin{split}
& \quad \mca(v_{1}^{n+1}+v_{2}^{n+1},u_{1}^{n+1})+\mca(v_{1}^{n+1}+v_{2}^{n+1},u_{2}^{n+1}) \\
& \geq \frac{1}{2\Delta t}  \sum_{i=1}^2 \norm{v_i^{n+1}}_{\mca}^2 - \frac{1}{2}\left ( \frac{1}{\Delta t} - \beta \right ) \sum_{i=1}^2 \left [ \norm{v_i^n}_{\mca}^2 - \norm{v_i^{n+1} - v_i^n}_{\mca}^2\right ] - \frac{1}{2\beta} \norm{u_2^{n+1} - u_2^n}_{\mca}^2. % - \frac{1}{2} \sum_{i=1}^2 \norm{u_i^{n+1}}_a^2. 
\end{split}
\label{eqn:proof-1}
\end{eqnarray}
Moreover, we have 
\begin{eqnarray}
\begin{split}
& \quad \left (\cfrac{u_{1}^{n+1}-u_{1}^{n}}{\Delta t}+\cfrac{u_{2}^{n+1}-u_{2}^{n}}{\Delta t},u_1^{n+1} \right )+
\left (\cfrac{u_{1}^{n+1}-u_{1}^{n}}{\Delta t}+\cfrac{u_{2}^{n+1}-u_{2}^{n}}{\Delta t},u_2^{n+1}\right )\\
=& \cfrac{1}{\Delta t}(u^{n+1}-u^n, u^{n+1})
= \cfrac{1}{2\Delta t} (\|u^{n+1}\|^2 - \|u^n\|^2 + \|u^{n+1} - u^n \|^2) \\
\geq &\cfrac{1}{2\Delta t} \left (\|u^{n+1}\|^2 - \|u^n\|^2 +  \sum_{i=1}^2\|u_i^{n+1} - u_i^n \|^2 - 2\gamma \|u_1^{n+1} - u_1^n \| \|u_2^{n+1} - u_2^n \| \right )\\
\geq &\cfrac{1}{2\Delta t} \left (\|u^{n+1}\|^2 - \|u^n\|^2 +  (1-\gamma) \sum_{i=1}^2\|u_i^{n+1} - u_i^n \|^2\right ). 
\end{split}
\label{eqn:proof-2}
\end{eqnarray}
Adding \eqref{eqn:proof-1} and \eqref{eqn:proof-2}, we obtain  
\begin{eqnarray*}
\begin{split}
0 & = \left (\cfrac{u_{1}^{n+1}-u_{1}^{n}}{\Delta t}+\cfrac{u_{2}^{n+1}-u_{2}^{n}}{\Delta t},u_1^{n+1} \right )+
\left (\cfrac{u_{1}^{n+1}-u_{1}^{n}}{\Delta t}+\cfrac{u_{2}^{n+1}-u_{2}^{n}}{\Delta t},u_2^{n+1}\right ) \\
& \quad + \mca(v_{1}^{n+1}+v_{2}^{n+1},u_{1}^{n+1})+\mca(v_{1}^{n+1}+v_{2}^{n+1},u_{2}^{n+1}) \\
& \geq \cfrac{1}{2\Delta t} \left (\|u^{n+1}\|^2 - \|u^n\|^2 +  \sum_{i=1}^2(1-\gamma) \|u_i^{n+1} - u_i^n \|^2\right) 
+ \frac{1}{2\Delta t} \sum_{i=1}^2 \norm{v_i^{n+1}}_{\mca}^2 \\
& \quad - \frac{1}{2}\left ( \frac{1}{\Delta t} - \beta \right ) \sum_{i=1}^2 \left [ \norm{v_i^n}_{\mca}^2 - \norm{v_i^{n+1} - v_i^n}_{\mca}^2\right ] - \frac{1}{2\beta} \norm{u_2^{n+1} - u_2^n}_{\mca}^2. % - \frac{1}{2} \sum_{i=1}^2 \norm{u_i^{n+1}}_a^2. 
\end{split}
\end{eqnarray*}
If the stability condition \eqref{eqn:pe-time-step} holds, then we obtain 
\begin{align*}
\tilde E^{n+1}(u,v) & =\norm{u^{n+1}}^2+ \sum_{i=1}^2 \norm{v_i^{n+1}}_{\mca}^2 \\
&\leq \|u^{n}\|^2 + (1-\beta \Delta t) \sum_{i=1}^2 \left [ \|v_i^{n}\|_{\mca}^{2}-\|v_i^{n+1}-v_i^{n}\|_{\mca}^{2}\right] -  (1-\gamma) \|u_1^{n+1} - u_1^n \|^2\\
& \quad - \left [ \underbrace{(1-\gamma) \|u_2^{n+1} - u_2^n \|^2 - \cfrac{\Delta t}{\beta}\|u_2^{n+1}-u_2^n\|^2_{\mca}}_{ \geq 0}\right ] \leq \|u^n\|^2 + \sum_{i=1}^2 \norm{v_i^{n}}_{\mca}^2 = \tilde E^n(u,v)
\end{align*}
for any $n \in \{0, 1,\cdots, N_T -1\}$. This completes the proof. 

\end{document}